# TESTING PREDICTOR CONTRIBUTIONS IN SUFFICIENT DIMENSION REDUCTION

By R. Dennis Cook[1]

*University of Minnesota*

We develop tests of the hypothesis of no effect for selected predictors in regression, without assuming a model for the conditional distribution of the response given the predictors. Predictor effects need not be limited to the mean function and smoothing is not required. The general approach is based on sufficient dimension reduction, the idea being to replace the predictor vector with a lower-dimensional version without loss of information on the regression. Methodology using sliced inverse regression is developed in detail.

**1. Introduction.** In full generality, the goal of a regression is to infer about the conditional distribution of the univariate response variable $Y$ given the $p \times 1$ vector of predictors $\mathbf{X}$: How does the conditional distribution of $Y|\mathbf{X}$ change with the value assumed by $\mathbf{X}$? Many different statistical contexts have been developed to address this issue. In this article we consider *sufficient dimension reduction* (SDR), the basic idea being to replace the predictor vector with its projection onto a subspace of the predictor space *without loss of information* on $Y|\mathbf{X}$. More formally, we seek subspaces $\mathcal{S}$ of the predictor space with the property that

$$Y \perp\!\!\!\perp \mathbf{X} | P_{\mathcal{S}}\mathbf{X}, \tag{1}$$

where $\perp\!\!\!\perp$ indicates independence, $P_{(\cdot)}$ stands for a projection operator in the standard inner product and, for future reference, $Q_{(\cdot)} = I_p - P_{(\cdot)}$. The statement is thus that $Y$ is independent of $\mathbf{X}$ given any value for $P_{\mathcal{S}}\mathbf{X}$. Subspaces with this property are called dimension reduction subspaces. Letting $k = \dim(\mathcal{S})$, a regression inquiry can then be limited to $k \le p$ new predictors, expressed as linear combinations of the original ones: $\mathbf{v}_1^T\mathbf{X}, \ldots, \mathbf{v}_k^T\mathbf{X}$, where

Received January 2003; revised May 2003.
[1]Supported in part by NSF Grant DMS-01-03983.

the basis $\{\mathbf{v}_1,\ldots,\mathbf{v}_k\}$ for $\mathcal{S}$ is often chosen so that the new predictors are uncorrelated.

When the intersection of all subspaces satisfying (1) also satisfies (1) it is called the *central subspace* (CS) [Cook (1994, 1996, 1998a)] and is denoted by $\mathcal{S}_{Y|\mathbf{X}}$. The central subspace, which is assumed to exist throughout this article, is a population metaparameter that can be taken as the parsimonious target of a dimension reduction inquiry. Its dimension $d = \dim(\mathcal{S}_{Y|\mathbf{X}})$ is called the *structural dimension* of the regression. There are several methods available that can be used to estimate the CS, including sliced inverse regression (SIR) [Li (1991)], sliced average variance estimation (SAVE) [Cook and Weisberg (1991)], graphical regression [Cook (1994, 1998a)], parametric inverse regression [Bura and Cook (2001b)] and partial SIR [Chiaromonte, Cook and Li (2002)] when categorical predictors are present. Cook and Weisberg (1999a) gave an introductory account of studying regressions via central subspaces.

Other dimension reduction methods estimate the *central mean subspace* [Cook and Li (2002)], which is a subspace of the CS that captures the mean function. These include ordinary least squares (OLS) and related methods based on convex objective functions, principal Hessian directions [Li (1992) and Cook (1998b)], iterative Hessian transformation [Cook and Li (2002)] and minimum average variance estimation [Xia, Tang, Li and Zhu (2002)]. In this article we are concerned only with the CS.

The estimation methods for the CS mentioned previously are all consistent under reasonable conditions when the dimension $d$ of the CS is known. Inference on $d$ is often based on hypothesis testing: Starting with $m = 0$, test the hypothesis $d = m$ versus $d > m$. If the test is rejected, increment $m$ by 1 and test again, stopping with the first nonsignificant result. This type of procedure is fairly common for estimating the dimension of a subspace [see, e.g., Rao (1965), page 472]. Once an estimate $\hat{d}$ is obtained, subsequent analysis, including choice of a first model, is typically guided by a summary plot of $Y$ versus the new predictors $\hat{\boldsymbol{\eta}}_1^T\mathbf{X},\ldots,\hat{\boldsymbol{\eta}}_{\hat{d}}^T\mathbf{X}$, where $\hat{\boldsymbol{\eta}}_j \in \mathbb{R}^p$ and $\{\hat{\boldsymbol{\eta}}_1,\ldots,\hat{\boldsymbol{\eta}}_{\hat{d}}\}$ is the estimated basis for $\mathcal{S}_{Y|\mathbf{X}}$. Examples of this process are available throughout the SDR literature. For recent examples, see Chen and Li (1998), Cook and Lee (1999) and Chiaromonte, Cook and Li (2002).

The ability to test the significance of subsets of predictors is often important in model-based regression, but is currently unavailable in SDR. In this article we develop tests of hypotheses involving statements of the form $P_{\mathcal{H}}\mathcal{S}_{Y|\mathbf{X}} = \mathcal{O}_p$, where $\mathcal{H}$ is a user-selected subspace of the predictor space that specifies the hypothesis, and $\mathcal{O}_p$ indicates the origin in $\mathbb{R}^p$. Partitioning $\mathbf{X}^T = (\mathbf{X}_1^T, \mathbf{X}_2^T)$, we imagine a typical application to test the hypothesis that $r$ selected predictors $\mathbf{X}_2$ do not contribute to the regression. Let the columns of the $p \times d$ matrix $\boldsymbol{\eta}$ be a basis for $\mathcal{S}_{Y|\mathbf{X}}$ and partition $\boldsymbol{\eta}^T = (\boldsymbol{\eta}_1^T, \boldsymbol{\eta}_2^T)$ according to the partition of $\mathbf{X}$. By definition of $\mathcal{S}_{Y|\mathbf{X}}$, $Y \perp\!\!\!\perp \mathbf{X}|\boldsymbol{\eta}^T\mathbf{X}$. We wish a test of the hypothesis $Y \perp\!\!\!\perp \mathbf{X}|\boldsymbol{\eta}_1^T\mathbf{X}_1$ so that the *coordinate subspace* $\text{Span}(\boldsymbol{\eta}_2)$ coincides with the origin, $\text{Span}(\boldsymbol{\eta}_2) = \mathcal{O}_r$. This can be expressed in terms of the statement $P_{\mathcal{H}}\mathcal{S}_{Y|\mathbf{X}} = \mathcal{O}_p$ by choosing $\mathcal{H} = \text{Span}((0, I_r)^T)$ to be the subspace of $\mathbb{R}^p$ corresponding to the coordinates $\mathbf{X}_2$ in question. Because we expect $\mathcal{H}$ will typically be chosen to target selected predictors, we refer to hypotheses of the form $P_{\mathcal{H}}\mathcal{S}_{Y|\mathbf{X}} = \mathcal{O}_p$ as *coordinate hypotheses*, although $\mathcal{H}$ need not correspond to a subset of predictors (coordinates). We let $r = \dim(\mathcal{H})$.

The following proposition gives a conditional independence interpretation of the statement $P_{\mathcal{H}}\mathcal{S}_{Y|\mathbf{X}} = \mathcal{O}_p$. Its proof is sketched in the Appendix.

PROPOSITION 1. *$P_{\mathcal{H}}\mathcal{S}_{Y|\mathbf{X}} = \mathcal{O}_p$ if and only if $Y \perp\!\!\!\perp P_{\mathcal{H}}\mathbf{X}|Q_{\mathcal{H}}\mathbf{X}$.*

Consequently, a coordinate hypothesis test can be regarded as a test of the hypothesis that, given $Q_{\mathcal{H}}\mathbf{X}$, the orthogonal part $P_{\mathcal{H}}\mathbf{X}$ of the predictor vector contains no information about the response. With $\mathcal{H} = \text{Span}((0, I_r)^T)$ the hypothesis $P_{\mathcal{H}}\mathcal{S}_{Y|\mathbf{X}} = \mathcal{O}_p$ is equivalent to the hypothesis that $Y$ and $\mathbf{X}_2$ are conditionally independent given $\mathbf{X}_1$, $Y \perp\!\!\!\perp \mathbf{X}_2|\mathbf{X}_1$.

In this article we consider three kinds of hypotheses that could be useful depending on the application-specific requirements:

1. *Marginal dimension hypotheses*—$d = m$ versus $d > m$;
2. *Marginal coordinate hypotheses*—$P_{\mathcal{H}}\mathcal{S}_{Y|\mathbf{X}} = \mathcal{O}_p$ versus $P_{\mathcal{H}}\mathcal{S}_{Y|\mathbf{X}} \neq \mathcal{O}_p$;
3. *Conditional coordinate hypotheses*—$P_{\mathcal{H}}\mathcal{S}_{Y|\mathbf{X}} = \mathcal{O}_p$ versus $P_{\mathcal{H}}\mathcal{S}_{Y|\mathbf{X}} \neq \mathcal{O}_p$ given $d$.

Marginal dimension hypotheses are considered extensively in the literature and are mentioned here for completeness. The other two forms are new and tests for them are developed in this article. Any of the dimension reduction methods mentioned previously (e.g., SIR, SAVE or PIR) could in principle be a foundation for tests of these hypotheses. In effect, graphical regression [Cook (1994, 1998a)] is built on our ability to assess coordinate hypotheses in a series of three-dimensional plots. In this article we use SIR to develop formal asymptotic tests of the two new hypotheses.

Our use of SIR to develop tests of hypotheses involving coordinate restrictions depends on rederiving it as the solution to a multivariate nonlinear least squares problem. This is done in Section 3.1 following further discussion of preliminary issues in Section 2. The population structure of SIR is related to the coordinate hypotheses in Section 3.2, and general results on test statistic construction are described in Section 4. In Sections 5 and 6 we develop the tests for the marginal and conditional coordinate hypotheses, including asymptotic null distributions and suggestions for implementation. Simulation results on level and power along with an illustrative data analysis are reported in Section 7. Concluding comments are given in Section 8, along with additional discussion of the literature and its relation to this work. To avoid interrupting the discussion, proofs for most results are given in the Appendix.

**2. Preparations.** We assume throughout this article that the data $(Y_i, \mathbf{X}_i)$, $i = 1, \ldots, n$, $\mathbf{X} \in \mathbb{R}^p$, are i.i.d. observations on $(Y, \mathbf{X})$, which has a joint distribution with finite fourth moments and $\mathbf{\Sigma} \equiv \text{Var}(\mathbf{X}) > 0$. In keeping with the usual SIR protocol, we assume also that the response has been discretized by constructing $h$ slices so that $Y$ takes values in $\{1, 2, \ldots, h\}$. The $j$th value of $Y$ is called the $j$th slice. This slicing step might be unnecessary if the response is naturally discrete or categorical.

Let the standardized predictors be denoted by

$$\mathbf{Z} = \mathbf{\Sigma}^{-1/2}(\mathbf{X} - \text{E}(\mathbf{X}))$$

with sample version

$$\hat{\mathbf{Z}}_{yj} = \hat{\mathbf{\Sigma}}^{-1/2}(\mathbf{X}_{yj} - \bar{\mathbf{X}}),$$

where subscript $(yj)$ indicates observation $j$ in slice $y$, $y = 1, \ldots, h$, $j = 1, \ldots, n_y$, $n = \sum_y n_y$, $\bar{\mathbf{X}} = \sum_{yj} \mathbf{X}_{yj}/n$ is the sample mean of the $\mathbf{X}_{yj}$'s,

$$\hat{\mathbf{\Sigma}} = \frac{1}{n} \sum_{y=1}^{h} \sum_{j=1}^{n_y} (\mathbf{X}_{yj} - \bar{\mathbf{X}})(\mathbf{X}_{yj} - \bar{\mathbf{X}})^T$$

is the usual sample covariance matrix and $\mathbf{\Sigma}^{-1/2}$ denotes the unique symmetric positive-definite square root of $\mathbf{\Sigma}^{-1}$. To allow use of the usual inner product in subsequent developments and without loss of generality, we work in the $\mathbf{Z}$ scale with central subspace $\mathcal{S}_{Y|\mathbf{Z}} = \mathbf{\Sigma}^{1/2} \mathcal{S}_{Y|\mathbf{X}}$ [Cook (1998a), Proposition 6.3], letting the columns of the $p \times d$ matrix $\boldsymbol{\gamma}$ be an orthonormal basis for $\mathcal{S}_{Y|\mathbf{Z}}$. Summations $\sum_{yj}$ with implicit limits $(yj)$ are always over $y = 1, \ldots, h$, $j = 1, \ldots, n_y$.

In practice coordinate hypotheses will typically be formulated in the original $\mathbf{X}$ scale by selecting an appropriate basis for $\mathcal{H}$. A coordinate hypothesis could then be stated as $\boldsymbol{\alpha}_x^T \boldsymbol{\eta} = 0$, where $\boldsymbol{\alpha}_x$ is the user-selected basis for $\mathcal{H}$ expressed as a $p \times r$ matrix of full column rank $r$, and $\boldsymbol{\eta}$ is a basis for $\mathcal{S}_{Y|\mathbf{X}}$. For example, to test if a selected subset of $r$ predictors contributes to the regression we can test if the rows of $\boldsymbol{\eta}$ corresponding to the $r$ predictors in question are all zero vectors. The matrix $\boldsymbol{\alpha}_x$ can then be chosen to select the appropriate rows of $\boldsymbol{\eta}$.

The hypothesis $\boldsymbol{\alpha}_x^T \boldsymbol{\eta} = 0$ holds if and only if $\boldsymbol{\alpha}^T (\mathbf{\Sigma}^{1/2} \boldsymbol{\eta}) = 0$, where $\boldsymbol{\alpha} = \mathbf{\Sigma}^{-1/2} \boldsymbol{\alpha}_x$ and the columns of $\mathbf{\Sigma}^{1/2} \boldsymbol{\eta}$ form a basis for $\mathcal{S}_{Y|\mathbf{Z}}$. A coordinate hypothesis in the $\mathbf{X}$ scale, $P_{\mathcal{H}} \mathcal{S}_{Y|\mathbf{X}} = \mathcal{O}_p$ with $\mathcal{H} = \text{Span}(\boldsymbol{\alpha}_x)$, can be restated in the $\mathbf{Z}$ scale as $P_{\mathcal{H}} \mathcal{S}_{Y|\mathbf{Z}} = \mathcal{O}_p$ with $\mathcal{H} = \text{Span}(\boldsymbol{\alpha})$. Thus by appropriate choice of basis, $\boldsymbol{\alpha}_x$ or $\boldsymbol{\alpha}$, we can work in either scale.

Back to the $\mathbf{Z}$ scale, without loss of generality we take the columns of

$$\boldsymbol{\alpha} = \boldsymbol{\Sigma}^{-1/2}\boldsymbol{\alpha}_x(\boldsymbol{\alpha}_x^T\boldsymbol{\Sigma}^{-1}\boldsymbol{\alpha}_x)^{-1/2} \tag{2}$$

to be an orthonormal basis for $\mathcal{H}$ in the remainder of this article. The hypothesis $P_{\mathcal{H}}\mathcal{S}_{Y|\mathbf{Z}} = \mathcal{O}_p$ holds if and only if $\mathcal{S}_{Y|\mathbf{Z}}$ is in the orthogonal complement of $\mathcal{H}$ and consequently under the hypothesis we must have $r \leq p - d$. Otherwise the hypothesis is certainly false.

## 3. SIR.

3.1. *Nonlinear least squares formulation.* The development of SIR as a means to estimate $\mathcal{S}_{Y|\mathbf{Z}}$ requires the following condition:

(C1) *Linearity condition*—$\mathrm{E}(\mathbf{Z}|P_{\mathcal{S}_{Y|\mathbf{Z}}}\mathbf{Z}) = P_{\mathcal{S}_{Y|\mathbf{Z}}}\mathbf{Z}$.

This condition, which is common in SDR, is equivalent to requiring that $\mathrm{E}(\mathbf{Z}|\boldsymbol{\gamma}^T\mathbf{Z})$ be a linear function of $\boldsymbol{\gamma}^T\mathbf{Z}$ [Cook (1998a), Proposition 4.2]. Li's (1991) design condition is equivalent to (C1), which applies to the marginal distribution of the predictors and not to the conditional distribution of $Y|\mathbf{Z}$ as is common in regression modeling. Consequently, we are free to use experimental design, one-to-one predictor transformations $\boldsymbol{\tau}$ or reweighting [Cook and Nachtsheim (1994)] to induce the condition when necessary without suffering complications when inferring about $Y|\mathbf{Z}$. Since we are not assuming a model for $Y|\mathbf{X}$, these adaptation methods need not change the fundamental issues in the regression. For example, because $Y|(\mathbf{X} = \mathbf{x})$ has the same distribution as $Y|(\boldsymbol{\tau}(\mathbf{X}) = \boldsymbol{\tau}(\mathbf{x}))$, predictor transformations just change the way in which the conditional distribution of $Y|\mathbf{X}$ is indexed. The linearity condition holds for elliptically contoured predictors. Additionally, Hall and Li (1993) showed that as $p$ increases with $d$ fixed the linearity condition holds to a reasonable approximation in many problems.

The linearity condition implies that the conditional means $\mathrm{E}(\mathbf{Z}|Y)$ lie in the CS for all values of $Y$ [Li (1991)]. We take this a step further and assume the following condition:

(C2) *Coverage condition*—$\mathrm{Span}\{\mathrm{E}(\mathbf{Z}|Y=y)|y=1,\ldots,h\} = \mathcal{S}_{Y|\mathbf{Z}}$,

so that the subspace spanned by the inverse conditional means coincides with the CS. This condition is also common in regression studies based on SIR. It requires in part that $h \geq d+1$. For subsequent tests on $d$ we require $h > d+1$.

For each value $y$ of $Y$ we can now find a vector $\boldsymbol{\rho}_y \in \mathbb{R}^d$ such that

$$\mathrm{E}(\mathbf{Z}|Y=y) = \boldsymbol{\gamma}\boldsymbol{\rho}_y,$$

where $\boldsymbol{\gamma}$ is the basis matrix for $\mathcal{S}_{Y|\mathbf{Z}}$ defined previously. Because $\mathrm{E}(\mathbf{Z}) = 0$ we must have $\mathrm{E}(\boldsymbol{\rho}_Y) = \sum_y f_y \boldsymbol{\rho}_y = 0$, where $f_y = \mathrm{Pr}(Y=y)$ is the probability of slice $y$. This suggests that for fixed $d$ estimates of $\boldsymbol{\gamma}$ and $\boldsymbol{\rho}_y$ can be constructed by minimizing the least squares loss function

$$L_d(\mathbf{B}, \mathbf{C}_y) \equiv \sum_{y=1}^{h} \sum_{j=1}^{n_y} \|\hat{\mathbf{Z}}_{yj} - \mathbf{B}\mathbf{C}_y\|^2$$

over $\mathbf{B}$ in the Stiefel manifold [Muirhead (1982), page 67] of all $p \times d$ semi-orthogonal matrices and over $\mathbf{C}_y \in \mathbb{R}^d$ subject to $\sum_y \hat{f}_y \mathbf{C}_y = 0$, where $\hat{f} = n_y/n$ is the observed fraction of observations falling in slice $y$. The values of $\mathbf{B}$ and $\mathbf{C}_y$ that minimize $L_d$ are then taken as the estimates $\hat{\boldsymbol{\gamma}}$ and $\hat{\boldsymbol{\rho}}_y$ of $\boldsymbol{\gamma}$ and $\boldsymbol{\rho}_y$, $y = 1\ldots,h$. Although we refer to $\hat{\boldsymbol{\gamma}}$ and $\hat{\boldsymbol{\rho}}_y$ as estimates, it may be, strictly speaking, more appropriate to think of them as solutions since they can be replaced by $\hat{\boldsymbol{\gamma}}\mathbf{H}^T$ and $\mathbf{H}\hat{\boldsymbol{\rho}}_y$, where $\mathbf{H}$ is any orthogonal matrix.

Minimizing $L_d$ results in the SIR estimate of $\mathcal{S}_{Y|\mathbf{Z}}$ when $d$ is regarded as known: Let $\bar{\mathbf{Z}}_y = \sum_j \hat{\mathbf{Z}}_{yj}/n_y$ be the average of the $\hat{\mathbf{Z}}_{yj}$ in slice $y$, and write $L_d$ as

$$L_d(\mathbf{B}, \mathbf{C}_y) = \sum_{yj} \|\hat{\mathbf{Z}}_{yj} - \bar{\mathbf{Z}}_y\|^2 + \sum_{yj} \|\bar{\mathbf{Z}}_y - \mathbf{B}\mathbf{C}_y\|^2.$$

For fixed $\mathbf{B}$ the minimum is attained by

$$\bar{\mathbf{C}}_y = \mathbf{B}^T \bar{\mathbf{Z}}_y, \qquad y = 1, \ldots, h.$$

Then minimizing $L_d(\mathbf{B}, \bar{\mathbf{C}}_y)$ over $\mathbf{B}$ yields the SIR estimate of $\mathcal{S}_{Y|\mathbf{Z}}$. To summarize the essential result, let $\widehat{\mathbf{M}} = \sum_y \hat{f}_y \bar{\mathbf{Z}}_y \bar{\mathbf{Z}}_y^T$ denote the sample covariance matrix of the slice means, and let $\hat{\lambda}_1 \geq \cdots \geq \hat{\lambda}_p$ denote the eigenvalues of $\widehat{\mathbf{M}}$. Then the columns of $\hat{\boldsymbol{\gamma}}$ are the eigenvectors corresponding to the first $d$ eigenvalues of $\widehat{\mathbf{M}}$, and $\hat{\boldsymbol{\rho}}_y = \hat{\boldsymbol{\gamma}}^T \bar{\mathbf{Z}}_y$, $y = 1, \ldots, h$.

The minimum value $\hat{L}_d \equiv L_d(\hat{\boldsymbol{\gamma}}, \hat{\boldsymbol{\rho}}_y)$, which we call the *residual sum of squares*, is

$$\hat{L}_d = \sum_{y=1}^{h} \sum_{j=1}^{n_h} \|\hat{\mathbf{Z}}_{yj} - \bar{\mathbf{Z}}_y\|^2 + n \sum_{j=d+1}^{p} \hat{\lambda}_j \tag{3}$$

for $d \leq p - 1$ and

$$\hat{L}_p = \sum_{y=1}^{h} \sum_{j=1}^{n_h} \|\hat{\mathbf{Z}}_{yj} - \bar{\mathbf{Z}}_y\|^2 \tag{4}$$

for $d = p$.

The usual SIR test statistic $T_n(m)$ for testing $d = m$ versus $d > m$, where $m < p$, can be found by comparing the residual sum of squares under the null hypothesis to that under the alternative,

$$T_n(m) = \hat{L}_m - \hat{L}_p = n \sum_{j=m+1}^{p} \hat{\lambda}_j. \tag{5}$$

Assuming that $\mathbf{X}$ has a multivariate normal distribution and implicitly assuming the coverage condition, Li (1991) proved that the distribution of $T_n(d)$ is asymptotically chi-squared with $(p-d)(h-d-1)$ degrees of freedom. Bura and Cook (2001a) proved that $T_n(d)$ has the same asymptotic distribution under the coverage and linearity conditions plus the following condition:

(C3) *Constant covariance condition*—$\mathrm{Var}(\mathbf{Z}|P_{\mathcal{S}_{Y|\mathbf{Z}}}\mathbf{Z}) = Q_{\mathcal{S}_{Y|\mathbf{Z}}}$,

where $Q_{\mathcal{S}_{Y|\mathbf{Z}}} = I_p - P_{\mathcal{S}_{Y|\mathbf{Z}}}$. This condition is equivalent to requiring that $\mathrm{Var}(\mathbf{Z}|P_{\mathcal{S}_{Y|\mathbf{Z}}}\mathbf{Z})$ be a nonrandom matrix. Normality of $\mathbf{X}$ implies the linearity and constant covariance conditions, but not the coverage condition. Bura and Cook (2001a) also proved that in general $T_n(d)$ is distributed as a weighted sum of independent chi-squared random variables and showed how to construct consistent estimates of the weights for use in practice.

In the next section we relate the coordinate hypothesis $P_{\mathcal{H}}\mathcal{S}_{Y|\mathbf{Z}} = \mathcal{O}_p$ to the population structure of SIR.

3.2. *Coordinate hypotheses and SIR.* Let $g_y = \sqrt{f_y}$, let $\boldsymbol{\mu}$ be the $p \times h$ matrix with columns $g_y \mathrm{E}(\mathbf{Z}|Y=y)$, $y = 1, \ldots, h$, and construct the singular value decomposition

$$\boldsymbol{\mu} = (\boldsymbol{\Gamma}_1 \ \boldsymbol{\Gamma}_0) \begin{pmatrix} \mathbf{D}_s & 0 \\ 0 & 0 \end{pmatrix} \begin{pmatrix} \boldsymbol{\Psi}_1^T \\ \boldsymbol{\Psi}_0^T \end{pmatrix}, \tag{6}$$

where $\boldsymbol{\Gamma} = (\boldsymbol{\Gamma}_1, \boldsymbol{\Gamma}_0)$ and $\boldsymbol{\Psi} = (\boldsymbol{\Psi}_1, \boldsymbol{\Psi}_0)$ are $p \times p$ and $h \times h$ orthogonal matrices, $\mathbf{D}_s$ is a $d \times d$ diagonal matrix of positive singular values and the various submatrices have the following dimensions:

$$\boldsymbol{\Gamma}_1 : p \times d, \qquad \boldsymbol{\Gamma}_0 : p \times p - d, \qquad \boldsymbol{\Psi}_1 : h \times d, \qquad \boldsymbol{\Psi}_0 : h \times h - d.$$

Under the linearity and coverage conditions, $\mathrm{Span}(\boldsymbol{\Gamma}_1) = \mathcal{S}_{Y|\mathbf{Z}}$ and so under these conditions we can take $\boldsymbol{\gamma} = \boldsymbol{\Gamma}_1$ as our basis for $\mathcal{S}_{Y|\mathbf{Z}}$.

The following two propositions relate coordinate hypotheses to the population structure of SIR. The proofs seem straightforward and are omitted.

PROPOSITION 2. *Assume that the linearity and coverage conditions hold. Then each of the following two conditions is equivalent to the coordinate hypothesis $P_{\mathcal{H}}\mathcal{S}_{Y|\mathbf{Z}} = \mathcal{O}_p$:*

(i) $Q_{\mathcal{H}}\boldsymbol{\Gamma}_1 = \boldsymbol{\Gamma}_1$.
(ii) $\mathcal{H} \subseteq \mathrm{Span}(\boldsymbol{\Gamma}_0)$.

*In addition, the coordinate hypothesis implies the following:*

(iii) $Q_{\mathcal{H}}\boldsymbol{\Gamma}_0 = \boldsymbol{\Gamma}_0(\boldsymbol{\Gamma}_0^T Q_{\mathcal{H}} \boldsymbol{\Gamma}_0)$.

(iv) $\mathbf{F}_{\mathcal{H}} \equiv \mathbf{\Gamma}_0^T Q_{\mathcal{H}} \mathbf{\Gamma}_0$ is a $(p-d) \times (p-d)$ symmetric idempotent matrix of rank $p-d-r$.

(v) $\mathbf{G}_{\mathcal{H}} \equiv I_{(p-d)} - \mathbf{F}_{\mathcal{H}} = \mathbf{\Gamma}_0^T P_{\mathcal{H}} \mathbf{\Gamma}_0$ is a $(p-d) \times (p-d)$ symmetric idempotent matrix of rank $r$.

PROPOSITION 3. *Assume that the linearity and coverage conditions hold. If $P_{\mathcal{H}} \mathcal{S}_{Y|\mathbf{Z}} = \mathcal{O}_p$, then the singular value decomposition of $\boldsymbol{\mu}$ is the same as that of $Q_{\mathcal{H}} \boldsymbol{\mu}$.*

**4. Test statistic construction.** In this section we discuss results that will facilitate construction of statistics for testing the two new hypotheses described in Section 1. Proceeding by analogy with the nonlinear least squares derivation of SIR described in Section 3.1, the test statistics will be constructed as the difference between the residual sums of squares under null and alternative hypotheses. The residual sum of squares under a dimension hypothesis $d = m$ can be written using (3)–(5) as

$$\hat{L}_m = \sum_{y=1}^{h} \sum_{j=1}^{n_y} \|\hat{\mathbf{Z}}_{yj} - \bar{\mathbf{Z}}_y\|^2 + T_n(m) \tag{7}$$

for $m = 0, \ldots, p$. Here we define $T_n(p) = 0$ so that (3) and (4) are both covered by (7).

We will also need the residual sum of squares under a coordinate constraint $P_{\mathcal{H}} \mathcal{S}_{Y|\mathbf{Z}} = \mathcal{O}_p$ and a dimension constraint $d = m$. Because $\boldsymbol{\Sigma}$ is typically unknown, it will have to be estimated for use in practice. Thus we let $\widehat{\mathcal{H}} = \mathrm{Span}(\hat{\boldsymbol{\Sigma}}^{-1/2} \boldsymbol{\alpha}_x)$.

To construct the residual sum of squares under coordinate and dimension constraints, write

$$L_m(\mathbf{B}, \mathbf{C}_y) = \sum_{yj} \|\hat{\mathbf{Z}}_{yj} - \bar{\mathbf{Z}}_y\|^2 + \sum_{yj} \|P_{\widehat{\mathcal{H}}}(\bar{\mathbf{Z}}_y - \mathbf{B}\mathbf{C}_y)\|^2 + \sum_{yj} \|Q_{\widehat{\mathcal{H}}}(\bar{\mathbf{Z}}_y - \mathbf{B}\mathbf{C}_y)\|^2.$$

Because $\mathbf{B}$ represents an orthonormal basis for $\mathcal{S}_{Y|\mathbf{Z}}$, we impose the constraint $P_{\widehat{\mathcal{H}}} \mathbf{B} = 0$, thus reducing $L_m$ to

$$L'_m(\mathbf{B}, \mathbf{C}_y) = \sum_{yj} \|\hat{\mathbf{Z}}_{yj} - \bar{\mathbf{Z}}_y\|^2 + \sum_{yj} \|P_{\widehat{\mathcal{H}}} \bar{\mathbf{Z}}_y\|^2 + \sum_{yj} \|Q_{\widehat{\mathcal{H}}}(\bar{\mathbf{Z}}_y - \mathbf{B}\mathbf{C}_y)\|^2, \tag{8}$$

where the prime on $L'_m$ indicates the imposition of the coordinate constraint. For fixed $\mathbf{B}$ with $\mathbf{B}^T \mathbf{B} = \mathbf{B}^T Q_{\widehat{\mathcal{H}}} \mathbf{B} = I_m$, the minimum is attained by $\bar{\mathbf{C}}_y =$ $\mathbf{B}^T Q_{\widehat{\mathcal{H}}} \bar{\mathbf{Z}}_y$, $y = 1, \ldots, h$. Consequently, with $m < p - r$,

$$\min_{(\mathbf{B}, \mathbf{C}_y)} \sum_{y=1}^{h} \sum_{j=1}^{n_y} \|Q_{\widehat{\mathcal{H}}}(\bar{\mathbf{Z}}_y - \mathbf{B}\mathbf{C}_y)\|^2 = \min_{\mathbf{B}} \sum_{yj} \|Q_{\widehat{\mathcal{H}}}(\bar{\mathbf{Z}}_y - \mathbf{B}\mathbf{B}^T Q_{\widehat{\mathcal{H}}} \bar{\mathbf{Z}}_y)\|^2$$

$$= n \sum_{j=m+1}^{p} \hat{\lambda}'_j = n \sum_{j=m+1}^{p-r} \hat{\lambda}'_j,$$

where $\hat{\lambda}'_1 \geq \cdots \geq \hat{\lambda}'_p$ are the eigenvalues of $Q_{\widehat{\mathcal{H}}} \widehat{\mathbf{M}} Q_{\widehat{\mathcal{H}}}$. The last equality follows since the last $r$ eigenvalues of $Q_{\widehat{\mathcal{H}}} \widehat{\mathbf{M}} Q_{\widehat{\mathcal{H}}}$ are all 0. If $m \geq p - r$, then

$$\min_{(\mathbf{B}, \mathbf{C}_y)} \sum_{yj} \|Q_{\widehat{\mathcal{H}}}(\bar{\mathbf{Z}}_y - \mathbf{B}\mathbf{C}_y)\|^2 = 0.$$

Substituting into $L'_m(\mathbf{B}, \mathbf{C}_y)$ given in (8) we obtain the residual sum of squares

$$\hat{L}'_m = \sum_{yj} \|\hat{\mathbf{Z}}_{yj} - \bar{\mathbf{Z}}_y\|^2 + \sum_{yj} \|P_{\widehat{\mathcal{H}}} \bar{\mathbf{Z}}_y\|^2 + T'_n(m), \quad (9)$$

where $T'_n(m) = n \sum_{j=m+1}^{p} \hat{\lambda}'_j$ and we adopt the convention that $T'_n(p) = 0$.

In the next two sections we use (7) and (9) to construct test statistics for the new hypotheses introduced in Section 1.

**5. Marginal coordinate hypotheses.** The marginal coordinate hypothesis $P_{\mathcal{H}} \mathcal{S}_{Y|\mathbf{Z}} = \mathcal{O}_p$ versus $P_{\mathcal{H}} \mathcal{S}_{Y|\mathbf{Z}} \neq \mathcal{O}_p$ can be used to test the contributions of selected predictors without requiring a statement concerning the dimension of $\mathcal{S}_{Y|\mathbf{Z}}$. The test statistic $T_n(\mathcal{H})$ is the difference between the residual sums of squares under the null and alternative hypotheses:

$$T_n(\mathcal{H}) = \hat{L}'_p - \hat{L}_p = n \operatorname{trace}(P_{\widehat{\mathcal{H}}} \widehat{\mathbf{M}} P_{\widehat{\mathcal{H}}}) \quad (10)$$

$$= \|\sqrt{n} \operatorname{vec}(\hat{\boldsymbol{\alpha}}^T \mathbb{Z}_n)\|^2, \quad (11)$$

where vec is the usual operator that maps a matrix into a vector by stacking its columns, $\hat{\boldsymbol{\alpha}} = \hat{\boldsymbol{\Sigma}}^{-1/2} \boldsymbol{\alpha}_x (\boldsymbol{\alpha}_x^T \hat{\boldsymbol{\Sigma}}^{-1} \boldsymbol{\alpha}_x)^{-1/2}$ is an orthonormal basis for $\widehat{\mathcal{H}}$ and $\mathbb{Z}_n$ is the $p \times h$ matrix with columns $\hat{g}_y \bar{\mathbf{Z}}_y$ so that $\widehat{\mathbf{M}} = \mathbb{Z}_n \mathbb{Z}_n^T$ and $\mathbb{Z}_n \xrightarrow{p} \boldsymbol{\mu}$. The representation of $T_n(\mathcal{H})$ given by (10) is what might be expected based on intuition: to test if $P_{\mathcal{H}} \mathcal{S}_{Y|\mathbf{Z}} = \mathcal{O}_p$ we consider the size of the projection of $\widehat{\mathbf{M}}$ onto the subspace specified by the hypothesis. Before using (11) to describe the asymptotic distribution of $T_n(\mathcal{H})$ we consider another form of the statistic that might provide additional insights.

Because $\operatorname{E}(\mathbf{Z}|Y) \in \mathcal{S}_{Y|\mathbf{Z}}$,

$$\boldsymbol{\nu}_y \equiv \boldsymbol{\Sigma}^{-1}(\operatorname{E}(\mathbf{X}|Y = y) - \operatorname{E}(\mathbf{X})) \in \mathcal{S}_{Y|\mathbf{X}}, \quad y = 1, \ldots, h.$$

Consequently, under the coordinate hypothesis we must have $\boldsymbol{\alpha}_x^T \boldsymbol{\nu}_y = 0$ for all $y$. Letting $\hat{\boldsymbol{\nu}}_y = \hat{\boldsymbol{\Sigma}}^{-1}(\bar{\mathbf{X}}_y - \bar{\mathbf{X}})$, the test statistic can be written in terms of the hypothesized estimates $\boldsymbol{\alpha}_x^T \hat{\boldsymbol{\nu}}_y$ of 0 as

$$T_n(\mathcal{H}) = \sum_{i=1}^{h} \hat{f}_y \hat{\boldsymbol{\nu}}_y^T \boldsymbol{\alpha}_x (\boldsymbol{\alpha}_x^T \hat{\boldsymbol{\Sigma}}^{-1} \boldsymbol{\alpha}_x^T)^{-1} \boldsymbol{\alpha}_x^T \hat{\boldsymbol{\nu}}_y.$$

5.1. *Asymptotic distributions.* A little setup is needed before we can describe the asymptotic distribution of $T_n(\mathcal{H})$. Define the indicator variable $J_y = 1$ if $Y$ is in slice $y$ and 0 otherwise, let $\boldsymbol{\beta}_y = \boldsymbol{\Sigma}^{-1} \operatorname{Cov}(\mathbf{X}, J_y)$ and let $\varepsilon_y = J_y - f_y - \boldsymbol{\beta}_y^T(\mathbf{X} - \operatorname{E}(\mathbf{X}))$ denote the population residual from the OLS fit of $J_y$ on $\mathbf{X}$. Let $\boldsymbol{\varepsilon}$ be the $h \times 1$ vector with elements $\varepsilon_y$, let $\mathbf{D}_g$ be the $h \times h$ diagonal matrix with $g_y$ on the diagonal and recall that $\boldsymbol{\alpha}$, the population version of $\hat{\boldsymbol{\alpha}}$, is defined by (2). Finally, let $\chi_1^2(D), \chi_2^2(D), \ldots, \chi_K^2(D)$ denote independent chi-squared random variables, where the degrees of freedom $D$ and $K$ vary with context.

THEOREM 1. *Assume that the linearity condition holds. Then, under the coordinate hypothesis $P_{\mathcal{H}} \mathcal{S}_{Y|\mathbf{Z}} = \mathcal{O}_p$, $\sqrt{n} \operatorname{vec}(\hat{\boldsymbol{\alpha}}^T \mathbb{Z}_n)$ converges in distribution to a normal random vector with mean 0 and covariance matrix*

(12) $$\boldsymbol{\Omega}_{\mathcal{H}} = \operatorname{E}(\mathbf{D}_g^{-1} \boldsymbol{\varepsilon} \boldsymbol{\varepsilon}^T \mathbf{D}_g^{-1} \otimes \boldsymbol{\alpha}^T \mathbf{Z} \mathbf{Z}^T \boldsymbol{\alpha}).$$

*Consequently, from* (11),

$$T_n(\mathcal{H}) \xrightarrow{\mathcal{L}} \sum_{i=1}^{hr} \omega_i \chi_i^2(1),$$

*where $\omega_1 \geq \omega_2 \geq \cdots \geq \omega_{hr}$ are the eigenvalues of $\boldsymbol{\Omega}_{\mathcal{H}}$.*

This theorem requires the linearity condition but not the coverage condition. If the coverage condition fails so SIR estimates a subspace $\mathcal{S}$ of $\mathcal{S}_{Y|\mathbf{Z}}$, it provides a test of $P_{\mathcal{H}} \mathcal{S} = \mathcal{O}_p$, but we will necessarily miss part of the CS. If the coverage condition holds, then SIR estimates the whole CS and the theorem provides a test of the complete hypothesis $P_{\mathcal{H}} \mathcal{S}_{Y|\mathbf{Z}} = \mathcal{O}_p$. As discussed later in Section 8, the test implied by this theorem might be useful even if the linearity condition fails.

If we have conditions C1–C3, then $\boldsymbol{\Omega}_{\mathcal{H}}$ can be simplified. Let $Q_g = I - \mathbf{g}\mathbf{g}^T$, where $\mathbf{g}$ denotes the $h \times 1$ vector with elements $g_y$.

COROLLARY 1. *If the linearity, coverage and constant covariance conditions hold, then*

(13) $$\boldsymbol{\Omega}_{\mathcal{H}} = (Q_g - \boldsymbol{\mu}^T \boldsymbol{\mu}) \otimes I_r$$

*and*

$$T_n(\mathcal{H}) \xrightarrow{\mathcal{L}} \sum_{j=1}^{h-1} \delta_j \chi_j^2(r),$$

*where $\delta_1 \geq \cdots \geq \delta_h = 0$ are the eigenvalues of $Q_g - \boldsymbol{\mu}^T \boldsymbol{\mu}$.*

5.2. *Implementation.* The test statistic $T_n(\mathcal{H})$ is the same for all versions of the test, but the reference distribution changes depending on conditions (C1)–(C3). In the most general case described in Theorem 1 we need to estimate the eigenvalues of the $hr \times hr$ covariance matrix $\boldsymbol{\Omega}_{\mathcal{H}}$ to construct the reference distribution. We can construct a consistent estimate $\widehat{\boldsymbol{\Omega}}_{\mathcal{H}}$ of $\boldsymbol{\Omega}_{\mathcal{H}}$ by substituting sample estimates for the unknown quantities:

(14) $$\widehat{\boldsymbol{\Omega}}_{\mathcal{H}} = \frac{1}{n} \sum_{y=1}^{h} \sum_{j=1}^{n_y} \mathbf{D}_{\hat{g}}^{-1} \hat{\boldsymbol{\varepsilon}}_{yj} \hat{\boldsymbol{\varepsilon}}_{yj}^T \mathbf{D}_{\hat{g}}^{-1} \otimes \hat{\boldsymbol{\alpha}}^T \hat{\mathbf{Z}}_{yj} \hat{\mathbf{Z}}_{yj}^T \hat{\boldsymbol{\alpha}},$$

where $\hat{\boldsymbol{\alpha}}$ and $\hat{\mathbf{Z}}_{yj}$ are as defined previously and $\mathbf{D}_{\hat{g}}$ is an $h \times h$ diagonal matrix with $\hat{g}_y$ on the diagonal. Also, $\hat{\boldsymbol{\varepsilon}}_{yj}$ is the $h \times 1$ vector of the residuals for observation $(yj)$, with one residual from each of the sample linear regressions of $J_y$ on $\mathbf{X}$. Letting $\hat{\omega}_i$ denote the eigenvalues of $\widehat{\boldsymbol{\Omega}}_{\mathcal{H}}$, a $p$-value for the coordinate hypothesis can be constructed by comparing the observed value of $T_n(\mathcal{H})$ to the percentage points of $\sum_{i=1}^{hr} \hat{\omega}_i \chi_i^2(1)$. There is a substantial literature on computing tail probabilities of the distribution of a linear combination of chi-squared random variables. See Field (1993) for an introduction. Alternatively, tail areas can usually be approximated adequately by using Satterthwaite's approximation.

We can proceed similarly under conditions (C1)–(C3). The $p$-value can be found by comparing $T_n(\mathcal{H})$ to the percentage points of the distribution of $\sum_{i=1}^{h-1} \hat{\delta}_i \chi_i^2(r)$, where $\hat{\delta}_1 \geq \cdots \geq \hat{\delta}_h = 0$ are the eigenvalues of

(15) $$\widetilde{\boldsymbol{\Omega}}_{\mathcal{H}} = (Q_{\hat{g}} - \mathbb{Z}_n^T \mathbb{Z}_n) \otimes I_r,$$

each with multiplicity $r$.

For ease of reference, we refer to the test using the weighted chi-squared reference distribution constructed from (14) as the *general* test. The test using reference distribution constructed from (15) will be called the *constrained* test. Both tests use the same statistic $T_n(\mathcal{H})$, but the reference distribution depends on applicable constraints, as given in Corollary 1.

**6. Conditional coordinate hypotheses.** The conditional coordinate hypothesis $P_{\mathcal{H}} \mathcal{S}_{Y|\mathbf{Z}} = \mathcal{O}_p$ versus $P_{\mathcal{H}} \mathcal{S}_{Y|\mathbf{Z}} \neq \mathcal{O}_p$ given $d$ might be useful when $d$ is specified as a modeling device, or when inference on $d$ using $T_n(m)$ results in a clear estimate. A test statistic $T_n(\mathcal{H}|d)$ can again be constructed as the

difference between the residual sum of squares under the null and alternative hypotheses:

$$T_n(\mathcal{H}|d) = \hat{L}'_d - \hat{L}_d$$

$$(16) \qquad = T_n(\mathcal{H}) - (T_n(d) - T'_n(d))$$

$$(17) \qquad = n\sum_{j=1}^{d}\hat{\lambda}_j - n\sum_{j=1}^{d}\hat{\lambda}'_j,$$

where $T'_n(d) = n\sum_{j=d+1}^{p}\hat{\lambda}'_j$ and $T_n(d) = n\sum_{j=d+1}^{p}\hat{\lambda}_j$ are as defined in (5) and (9). Form (17) gives one way to compute the statistic and shows that it depends on the largest $d$ eigenvalues of $Q_{\widehat{\mathcal{H}}}\widehat{\mathbf{M}}Q_{\widehat{\mathcal{H}}}$ and $\widehat{\mathbf{M}}$ for the null and alternative hypotheses. In contrast, the usual SIR statistic $T_n(m)$ depends on the smallest $p-m$ eigenvalues of $\widehat{\mathbf{M}}$. Form (16) will be easier to work with when developing the asymptotic distribution of $T_n(\mathcal{H}|d)$ because it allows us to use some known results. To develop the asymptotic distribution of $T_n(\mathcal{H}|d)$ we consider first the asymptotic distributions of $T'_n(d)$ and $T_n(d) - T'_n(d)$ because these are components of $T_n(\mathcal{H}|d)$ and may be of interest in their own right. For instance, $T'_n(m) = \hat{L}'_m - \hat{L}'_p$ and thus it can be viewed as a test statistic for a dimension hypothesis given a coordinate constraint.

6.1. *Asymptotic distribution of $T'_n(d)$.* The asymptotic distribution of $T'_n(d)$ can be found by using results of Bura and Cook (2001a). Define

$$\sqrt{n}\mathbf{U}_n \equiv \sqrt{n}\mathbf{\Gamma}_0^T(\mathbb{Z}_n - \boldsymbol{\mu})\boldsymbol{\Psi}_0 = \sqrt{n}\mathbf{\Gamma}_0^T\mathbb{Z}_n\boldsymbol{\Psi}_0.$$

Bura and Cook [(2001a), equations (8)–(13) and associated discussion] first used the general results of Eaton and Tyler (1994) on the asymptotic distribution of singular values of a random matrix to conclude that the asymptotic distribution of $T_n(d)$ is the same as that of $n\|\mathbf{U}_n\|^2$. They then established that

$$\sqrt{n}\,\text{vec}(\mathbb{Z}_n - \boldsymbol{\mu}) \xrightarrow{\mathcal{L}} N_{ph}(0, \boldsymbol{\Delta})$$

and thus that

$$(18) \qquad \sqrt{n}\,\text{vec}(\mathbf{U}_n) \xrightarrow{\mathcal{L}} N_{(p-d)(h-d)}(0, (\boldsymbol{\Psi}_0^T \otimes \boldsymbol{\Gamma}_0^T)\boldsymbol{\Delta}(\boldsymbol{\Psi}_0 \otimes \boldsymbol{\Gamma}_0)),$$

where the $hp \times hp$ matrix $\boldsymbol{\Delta}$ is as defined by Bura and Cook (2001a). It can be represented as an $h \times h$ array of $p \times p$ matrices $\boldsymbol{\Delta}_{ss} = I_p f_s + (1 - 2f_s)\boldsymbol{\Sigma}_{\mathbf{Z}|s}$ and $\boldsymbol{\Delta}_{ts} = g_t g_s(I_p - \boldsymbol{\Sigma}_{\mathbf{Z}|t} - \boldsymbol{\Sigma}_{\mathbf{Z}|s})$, where $\boldsymbol{\Sigma}_{\mathbf{Z}|s} = \text{Var}(\mathbf{Z}|Y=s)$, $s,t = 1,\ldots,h$. Thus

$$T_n(d) \xrightarrow{\mathcal{L}} \sum_{i=1}^{(p-d)(h-d)} \omega_i \chi_i^2(1),$$

where $\omega_1 \geq \omega_2 \geq \cdots \geq \omega_{(p-d)(h-d)}$ are the eigenvalues of the covariance matrix in the asymptotic distribution of $\sqrt{n}\,\text{vec}(\mathbf{U}_n)$ given in (18).

The asymptotic distribution of $T'_n(d)$ can be found similarly. Define
$$\sqrt{n}\mathbf{U}'_n = \sqrt{n}\boldsymbol{\Gamma}_0^T(Q_{\widehat{\mathcal{H}}}\mathbb{Z}_n - Q_{\mathcal{H}}\boldsymbol{\mu})\boldsymbol{\Psi}_0$$
$$= \sqrt{n}\boldsymbol{\Gamma}_0^T(Q_{\widehat{\mathcal{H}}}\mathbb{Z}_n)\boldsymbol{\Psi}_0,$$
where the second equality follows because $\boldsymbol{\mu}\boldsymbol{\Psi}_0 = 0$ from the singular value decomposition (6). It follows from Eaton and Tyler (1994) that $T'_n(d)$ and $n\|\text{vec}(\mathbf{U}'_n)\|^2$ are asymptotically equivalent because, from Proposition 3, $\boldsymbol{\mu}$ and $Q_{\mathcal{H}}\boldsymbol{\mu}$ have the same singular value decomposition. Now,
$$\sqrt{n}\,\text{vec}(\mathbf{U}'_n) = (I_{h-d} \otimes \boldsymbol{\Gamma}_0^T Q_{\widehat{\mathcal{H}}})\sqrt{n}\,\text{vec}(\mathbb{Z}_n\boldsymbol{\Psi}_0).$$
Since $\text{vec}(\mathbb{Z}_n\boldsymbol{\Psi}_0) \xrightarrow{p} 0$, it follows that $\sqrt{n}\,\text{vec}(\mathbb{Z}_n\boldsymbol{\Psi}_0)$ converges in distribution. Because $Q_{\widehat{\mathcal{H}}}$ converges in probability to $Q_{\mathcal{H}}$, it follows from Slutsky's theorem that we can replace $\widehat{\mathcal{H}}$ with $\mathcal{H}$ in $\sqrt{n}\,\text{vec}(\mathbf{U}'_n)$ without affecting its asymptotic distribution. Consequently, $\sqrt{n}\,\text{vec}(\mathbf{U}'_n)$ is asymptotically equivalent to
$$(I_{h-d} \otimes \boldsymbol{\Gamma}_0^T Q_{\mathcal{H}})\sqrt{n}\,\text{vec}(\mathbb{Z}_n\boldsymbol{\Psi}_0) = \sqrt{n}\,\text{vec}(\boldsymbol{\Gamma}_0^T Q_{\mathcal{H}}\mathbb{Z}_n\boldsymbol{\Psi}_0)$$
$$= \sqrt{n}\,\text{vec}(\mathbf{F}_{\mathcal{H}}\boldsymbol{\Gamma}_0^T \mathbb{Z}_n\boldsymbol{\Psi}_0)$$
$$= (I_{h-d} \otimes \mathbf{F}_{\mathcal{H}})\sqrt{n}\,\text{vec}(\mathbf{U}_n),$$
where the second equality follows from parts (iii) and (iv) of Proposition 2. Consequently, the asymptotic distribution of $T'_n(d)$ is the same as that of $n\|(I_{h-d} \otimes \mathbf{F}_{\mathcal{H}})\text{vec}(\mathbf{U}_n)\|^2$, which can be determined from the asymptotic distribution of $\sqrt{n}\,\text{vec}(\mathbf{U}_n)$ given in (18). This enables us to conclude the following.

PROPOSITION 4.
$$T'_n(d) \xrightarrow{\mathcal{L}} \sum_{i=1}^{(p-d)(h-d)} \omega_i \chi_i^2(1),$$
where $\omega_1 \geq \omega_2 \geq \cdots \geq \omega_{(p-d)(h-d)}$ are the eigenvalues of
$$\boldsymbol{\Omega}'_d = (\boldsymbol{\Psi}_0^T \otimes \mathbf{F}_{\mathcal{H}}\boldsymbol{\Gamma}_0^T)\boldsymbol{\Delta}(\boldsymbol{\Psi}_0 \otimes \boldsymbol{\Gamma}_0\mathbf{F}_{\mathcal{H}}).$$

Additionally, the following corollary follows from Bura and Cook [(2001a), Theorem 2] and the fact that $\mathbf{F}_{\mathcal{H}}$ is a symmetric idempotent matrix of rank $p - d - r$ [Proposition 2(iv)].

COROLLARY 2. *Assume that the linearity, coverage and constant covariance conditions hold. Then $T'_n(d)$ is distributed asymptotically as a chi-squared random variable with $(p-d-r)(h-d-1)$ degrees of freedom.*

Given that $P_{\mathcal{H}}\mathcal{S}_{Y|\mathbf{Z}} = \mathcal{O}_p$, $\boldsymbol{\mu}\boldsymbol{\mu}^T$ and $Q_{\mathcal{H}}\boldsymbol{\mu}\boldsymbol{\mu}^T Q_{\mathcal{H}}$ have the same rank $d$. Consequently, we might expect $T_n(d) - T'_n(d)$ to reflect little more than random variation. Consider the orthogonal decomposition

$$n\|\operatorname{vec}(\mathbf{U}_n)\|^2 = n\|(I_{h-d} \otimes \mathbf{F}_{\mathcal{H}})\operatorname{vec}(\mathbf{U}_n)\|^2 + n\|(I_{h-d} \otimes \mathbf{G}_{\mathcal{H}})\operatorname{vec}(\mathbf{U}_n)\|^2,$$

where $\mathbf{G}_{\mathcal{H}}$ is as defined in Proposition 2(v). As discussed previously in this section, the left-hand side is asymptotically equivalent to $T_n(d)$ and the first term on the right-hand side is asymptotically equivalent to $T'_n(d)$. Thus, the second term is asymptotically equivalent to $T_n(d) - T'_n(d)$:

(19) $$T_n(d) - T'_n(d) = n\|(I_{h-d} \otimes \mathbf{G}_{\mathcal{H}})\operatorname{vec}(\mathbf{U}_n)\|^2 + o_p(1).$$

The next corollary gives the asymptotic distribution of $T_n(d) - T'_n(d)$ under conditions C1–C3. Its proof parallels that of Corollary 2 and is omitted.

COROLLARY 3. *Assume that the linearity, coverage and constant covariance conditions hold. Then $T_n(d) - T'_n(d)$ is distributed asymptotically as a chi-squared random variable with $r(h-d-1)$ degrees of freedom.*

6.2. *Asymptotic distribution of $T_n(\mathcal{H}|d)$.* The asymptotic distribution of $T_n(\mathcal{H}|d)$ can be found under the coordinate hypothesis by using the following proposition. The proof given in the Appendix relies on (19).

PROPOSITION 5.
$$T_n(\mathcal{H}|d) = \|(\boldsymbol{\Psi}_1^T \otimes I_r)\sqrt{n}\operatorname{vec}(\hat{\boldsymbol{\alpha}}^T \mathbb{Z}_n)\|^2 + o_p(1).$$

Using this proposition in combination with Theorem 1 gives the following theorem.

THEOREM 2. *Assume that the linearity and coverage conditions hold. Then, under the coordinate hypothesis $P_{\mathcal{H}}\mathcal{S}_{Y|\mathbf{Z}} = \mathcal{O}_p$, $(\boldsymbol{\Psi}_1^T \otimes I_r)\sqrt{n}\operatorname{vec}(\hat{\boldsymbol{\alpha}}^T \mathbb{Z}_n)$ converges in distribution to a normal random vector with mean $0$ and covariance matrix*

$$\boldsymbol{\Omega}_{\mathcal{H}|d} = \mathrm{E}(\boldsymbol{\Psi}_1^T \mathbf{D}_g^{-1}\boldsymbol{\varepsilon}\boldsymbol{\varepsilon}^T \mathbf{D}_g^{-1}\boldsymbol{\Psi}_1 \otimes \boldsymbol{\alpha}^T \mathbf{Z}\mathbf{Z}^T \boldsymbol{\alpha}).$$

*Consequently,*

$$T_n(\mathcal{H}|d) \overset{\mathcal{L}}{\to} \sum_{i=1}^{dr} \omega_i \chi_i^2(1),$$

*where $\omega_1 \geq \omega_2 \geq \cdots \geq \omega_{dr}$ are the eigenvalues of $\boldsymbol{\Omega}_{\mathcal{H}|d}$.*

It may be useful when reading this theorem to recall that $r \leq p - d$ for a meaningful coordinate hypothesis. In particular, $\mathbf{\Omega}_{\mathcal{H}|d}$ is not defined when $d = p$.

As in Section 5, if conditions (C1)–(C3) hold, then $\mathbf{\Omega}_{\mathcal{H}|d}$ can be simplified:

COROLLARY 4. *If the linearity, coverage and constant covariance conditions hold then*

$$\mathbf{\Omega}_{\mathcal{H}|d} = (I_d - \mathbf{D}_\lambda) \otimes I_r$$

*and*

$$T_n(\mathcal{H}|d) \xrightarrow{\mathcal{L}} \sum_{j=1}^{d}(1 - \lambda_j)\chi_j^2(r),$$

*where $\lambda_1 \geq \cdots \geq \lambda_d > 0$ are the nonzero eigenvalues of $\boldsymbol{\mu}\boldsymbol{\mu}^T$ and $\mathbf{D}_\lambda$ is a diagonal matrix with diagonal elements $\lambda_j$, $j = 1, \ldots, d$.*

The generalized inverse of $\mathbf{\Omega}_{\mathcal{H}|d}$ in Corollary 4 could be used to construct a Wald test statistic with an asymptotic chi-squared distribution under the coordinate hypothesis of Theorem 2. A similar comment applies to (13) under the coordinate hypothesis of Theorem 1.

6.3. *Implementation.* The results of Theorem 2 can be implemented in a manner similar to the implementation of Theorem 1 described in Section 5.2. A consistent estimate $\hat{\mathbf{\Psi}}_1$ of $\mathbf{\Psi}_1$ can be constructed from the singular value decomposition of $\mathbb{Z}_n$ just as $\mathbf{\Psi}_1$ is obtained from the singular value decomposition of $\boldsymbol{\mu}$ given in (6). A consistent estimate of $\mathbf{\Omega}_{\mathcal{H}|d}$ can then be constructed as

$$(20) \qquad \widehat{\mathbf{\Omega}}_{\mathcal{H}|d} = (\hat{\mathbf{\Psi}}_1^T \otimes I_r)\widehat{\mathbf{\Omega}}_\mathcal{H}(\hat{\mathbf{\Psi}}_1 \otimes I_r),$$

where $\widehat{\mathbf{\Omega}}_\mathcal{H}$ is as given in (14). Similarly, the asymptotic reference distribution of Corollary 4 can be estimated by substituting the largest $d$ eigenvalues $\hat{\lambda}_1, \ldots, \hat{\lambda}_d$ of $\widehat{\mathbf{M}}$ for $\lambda_1, \ldots, \lambda_d$, which amounts to estimating $\mathbf{\Omega}_{\mathcal{H}|d}$ by using

$$(21) \qquad \widetilde{\mathbf{\Omega}}_{\mathcal{H}|d} = (I_d - \mathbf{D}_{\hat{\lambda}}) \otimes I_r,$$

where $\mathbf{D}_{\hat{\lambda}}$ is a diagonal matrix with diagonal elements $\hat{\lambda}_j$, $j = 1, \ldots, d$.

Following the terminology for tests of marginal coordinate hypotheses, we refer to the test using reference distribution constructed from (20) as the *general* test. The *constrained* test uses the weighted chi-squared reference distribution based on (21). These two tests use the same statistic $T_n(\mathcal{H}|d)$; only the reference distribution changes.

**7. Simulation results and data analysis.** Simulation studies were conducted to insure that the asymptotic tests behave as expected and to provide a little insight about their operating characteristics. Each study was based on one of the following two models:

$$Y = X_1 + \varepsilon, \tag{22}$$

$$Y = \frac{X_1}{0.5 + (X_2 + 1.5)^2} + \delta. \tag{23}$$

The number of observations $n$, the number and distribution of the predictors $\mathbf{X}$ and the distributions of the errors $\varepsilon$ and $\delta$ depend on the simulation. To avoid inadvertent tuning by choice of the number of slices, every simulation run used $h = 5$ slices. Test results were tabulated over 1000 replications for each sampling configuration.

7.1. *Estimated versus nominal levels.* In this section we report some representative results to compare estimated and nominal levels. The estimates were obtained by counting the number of $p$-values that were less than or equal to a nominal level in the 1000 replications for each sample configuration. These $p$-values were obtained by applying the tests to a predictor not represented in the mean function of the model, so $r = 1$.

Estimated levels of all seven statistics described here are shown in Table 1 for simulations from model (22) with five i.i.d. standard normal predictors, an independent normal error and various sample sizes. For instance, the estimated levels shown in subtable A are for the test statistic $T_n(\mathcal{H})$ with its general reference distribution. The $p$-values were computed by comparing $T_n(\mathcal{H})$ to the quantiles of the weighted chi-squared distribution constructed by using the covariance matrix in (14). The results seem quite good for $n = 100$ and 200. Tests based on an estimated weighted chi-squared distribution (subtables A–D) tend to be liberal. This conclusion held up throughout all the simulations of test level conducted. The performance of the chi-squared statistics (subtables E–G), which tended to be conservative, was similar to that reported by Bura and Cook (2001a) for $T_n(d)$. The statistics $T'_n(d)$ and $T_n(d) - T'_n(d)$ were included in Table 1 to provide numerical support for the asymptotic calculations described previously. An investigation of possible roles for them in data analysis is outside the scope of this report. Discussion in the remainder of this section is confined to tests of the marginal and conditional coordinate hypotheses.

A substantial increase in the number of predictors typically required that the sample size be increased to achieve consistent agreement between the estimated and nominal levels. Shown in Table 2 are estimated levels for the two general and two constrained tests based on model (23) with $p = 10$ independent standard normal predictors. The agreement between the estimated and nominal levels for $n = 400$ and 800 seems quite good. Comparing the results for $T_n(\mathcal{H})$ with those for $T_n(\mathcal{H}|d)$ at $n = 50, 100$ suggests that tests based on $T_n(\mathcal{H}|d)$ need somewhat larger sample sizes to achieve similar agreement. This might be because use of $T_n(\mathcal{H}|d)$ requires an estimate of $\mathbf{\Psi}_1$ that is not required to use $T_n(\mathcal{H})$ [see (20)].

The two general tests, one for marginal coordinate hypotheses and one for conditional coordinate hypotheses, will probably be the most useful in practice since they require the fewest assumptions. In comparison, the cor-responding constrained tests achieved similar agreement between the estimated and nominal levels with somewhat smaller sample sizes.

TABLE 1
*Estimated level of nominal $1, 5, 10$ and $15\%$ tests based on various statistics and reference distributions for model (22) with $p = 5$ independent standard normal predictors and $\varepsilon = 0.2N(0, 1)$*

|     | Nominal level (%) | | | |
| --- | --- | --- | --- | --- |
| $n$ | 1 | 5 | 10 | 15 |
| (A) $T_n(\mathcal{H})$ with $\widehat{\mathbf{\Omega}}_\mathcal{H}$ (14) | | | | |
| 50  | 2.8 | 9.1 | 16.9 | 21.7 |
| 100 | 1.1 | 6.1 | 11.5 | 18.8 |
| 200 | 1.0 | 5.3 | 11.4 | 16.7 |
| (B) $T_n(\mathcal{H})$ with $\widetilde{\mathbf{\Omega}}_\mathcal{H}$ (15) | | | | |
| 50  | 2.1 | 8.1 | 15.2 | 20.1 |
| 100 | 1.0 | 5.6 | 10.9 | 17 |
| 200 | 0.9 | 5.3 | 10.3 | 16.3 |
| (C) $T_n(\mathcal{H}|d)$ with $\widehat{\mathbf{\Omega}}_{\mathcal{H}|d}$ (20) | | | | |
| 50  | 4.2 | 10.2 | 16.4 | 23 |
| 100 | 2.4 | 7.3 | 12.3 | 18.5 |
| 200 | 1.7 | 5.3 | 10.4 | 14.9 |
| (D) $T_n(\mathcal{H}|d)$ with $\widetilde{\mathbf{\Omega}}_{H|d}$ (21) | | | | |
| 50  | 3.0 | 8.5 | 14.7 | 20.6 |
| 100 | 1.9 | 6.3 | 12.2 | 17.6 |
| 200 | 1.6 | 5.2 | 9.9 | 14.6 |
| (E) $T_n(d) \sim \chi^2(12)$ | | | | |
| 50  | 0.4 | 4.6 | 11 | 17.2 |
| 100 | 0.9 | 4.1 | 9.1 | 14.7 |
| 200 | 1.4 | 4.9 | 9.7 | 14.1 |
| (F) $T'_n(d) \sim \chi^2(9)$ | | | | |
| 50  | 0.5 | 4.8 | 10.3 | 15.4 |
| 100 | 0.7 | 4.2 | 9.2 | 14.8 |
| 200 | 1.0 | 4.8 | 9.3 | 15.1 |
| (G) $T_n(d) - T'_n(d) \sim \chi^2(3)$ | | | | |
| 50  | 1.5 | 5.3 | 12.2 | 17.9 |
| 100 | 0.9 | 4.5 | 9.1 | 15.2 |
| 200 | 0.9 | 4.9 | 10.0 | 16.0 |

TABLE 2
*Estimated levels from model* (23) *with*
$p = 10$ *independent* $N(0,1)$ *predictors and*
$\delta = 0.2 N(0,1)$

|  | Nominal level (%) | | | |
|---|---|---|---|---|
| $n$ | 1 | 5 | 10 | 15 |
| (A) $T_n(\mathcal{H})$ with $\widehat{\boldsymbol{\Omega}}_{\mathcal{H}}$ (14) | | | | |
| 50 | 3.3 | 11.6 | 22.8 | 31.9 |
| 100 | 1.8 | 7.8 | 16.0 | 21.1 |
| 200 | 2.2 | 7.0 | 13.0 | 18.1 |
| 400 | 1.3 | 4.8 | 9.8 | 15.1 |
| 800 | 1.4 | 5.8 | 10.3 | 14.9 |
| (B) $T_n(\mathcal{H})$ with $\widetilde{\boldsymbol{\Omega}}_{\mathcal{H}}$ (15) | | | | |
| 50 | 2.9 | 9.8 | 19.2 | 29.5 |
| 100 | 1.3 | 7.2 | 14.2 | 19.9 |
| 200 | 1.9 | 6.9 | 12.2 | 17.7 |
| 400 | 1.2 | 4.8 | 10.0 | 14.4 |
| 800 | 1.4 | 5.9 | 10.1 | 14.8 |
| (C) $T_n(\mathcal{H}|d)$ with $\widehat{\boldsymbol{\Omega}}_{\mathcal{H}|d}$ (20) | | | | |
| 50 | 7.2 | 17.7 | 26.3 | 31.1 |
| 100 | 4.1 | 9.2 | 15.5 | 20.8 |
| 200 | 2.0 | 8.0 | 14.4 | 19.7 |
| 400 | 0.8 | 5.1 | 10.5 | 15.0 |
| 800 | 0.8 | 4.6 | 10.5 | 14.5 |
| (D) $T_n(\mathcal{H}|d)$ with $\widetilde{\boldsymbol{\Omega}}_{\mathcal{H}|d}$ (21) | | | | |
| 50 | 5.7 | 15.3 | 23.9 | 30.2 |
| 100 | 3.2 | 8.2 | 14.9 | 20.0 |
| 200 | 1.6 | 7.6 | 14.0 | 19.1 |
| 400 | 0.9 | 4.4 | 10.2 | 14.6 |
| 800 | 0.8 | 4.9 | 10.4 | 14.5 |

The results in Table 3 are intended to give some idea about the impact of the predictor distribution on the actual level of the two general tests. The subtables are designated as A and C to correspond to their designations in Tables 1 and 2. The simulation setup leading to Table 3 was repeated with other predictor distributions, including the $t$ distribution with five degrees of freedom and the uniform $(-2, 2)$ distribution. The results for these predictor distributions were quite similar to the results in Table 3.

Over the range of simulations represented in this study it was observed that the estimated level of a nominal 1% test was nearly always between 1 and 5% and the estimated level of a nominal 5% test was nearly always between 5 and 10%. No simulations were conducted with more than 12 predictors or more than 800 observations.

TABLE 3
*Estimated test levels from model* (23)
*with* 10 *independent* $\chi^2(4)$ *predictors and*
$\delta = 0.2N(0,1)$

|  | Nominal level (%) | | | |
|---|---|---|---|---|
| $n$ | 1 | 5 | 10 | 15 |
| (A) $T_n(\mathcal{H})$ with $\widehat{\boldsymbol{\Omega}}_{\mathcal{H}}$ (14) | | | | |
| 100 | 2.5 | 8.2 | 13.4 | 19.3 |
| 200 | 1.1 | 5.4 | 11.0 | 15.7 |
| 400 | 1.2 | 6.3 | 11.5 | 17.2 |
| 800 | 1.3 | 5.3 | 10.4 | 15.4 |
| (C) $T_n(\mathcal{H}|d)$ with $\widehat{\boldsymbol{\Omega}}_{\mathcal{H}|d}$ (20) | | | | |
| 100 | 1.6 | 7.2 | 12.7 | 18.7 |
| 200 | 1.6 | 7.3 | 12.5 | 18.9 |
| 400 | 0.7 | 3.2 | 7.9 | 12.9 |
| 800 | 1.0 | 5.6 | 10.0 | 16.0 |

7.2. *Power.* In this section we report results from a power study to gain insight into the operating characteristics of the proposed tests. It is not difficult to find examples where the power is near 1, the nominal level or anywhere between these extremes. To provide a benchmark for interpretation, the standard linear model $t$-test was included in the study.

The results reported in Table 4 are from model (22) with five independent standard normal predictors, $n = 200$ and three different errors $\varepsilon$. For each model configuration, the power of the standard $t$-test for the hypothesis that the coefficient of $X_1$ equals 0, and the power of the general marginal coordinate test for $X_1$, were estimated by computing the fraction of rejections in 1000 replications. The first column of Table 4 indicates the test. The second column indicates the nature of the error and will be described shortly. The third and fourth columns give the estimated power ($PR$) at the nominal 1 and 5% levels. The differences between the estimated and nominal levels for all tests in Table 4 were found to be roughly as those of Table 1.

To provide some information about estimation in addition to that for testing, we also computed the absolute sample correlations $c$ between $X_1$ and the fitted values from the OLS fit of $Y$ on $\mathbf{X}$, including an intercept, and between $X_1$ and the first SIR predictor. The 0.05, 0.5 and 0.95 quantiles $c_{0.05}$, $c_{0.5}$ and $c_{0.9}$ of the empirical distributions of these absolute correlations are given in columns 5–7 of Table 4.

*Table* 4(A), $\varepsilon = \sigma N(0,1)$. For $\sigma \leq 2$ the two procedures were observed to yield essentially identical results. Both tests rejected in all 1000 replications, and the absolute correlations were all quite high. The results for $\sigma = 1$ are shown in the first two rows. The $t$-test was observed to be the clear winner for $\sigma \geq 3$; the results for $\sigma = 6.4$ are shown in the third and fourth rows of this table. The qualitative nature of these results should perhaps not be surprising since the $t$-test has the home field advantage with a homoscedastic normal error. The estimated powers at 1 and 5% of the general conditional coordinate test $T(\mathcal{H}|d=1)$ were observed to be 0.275 and 0.469 for $\sigma = 6.4$. Comparing these results with the corresponding results in the table suggests that a substantial part of the power differences between the $t$- and $T_n(\mathcal{H})$-tests can be attributed to the differential information on dimension.

*Table* 4(B), $\varepsilon = 6.4(\chi^2(D) - D)/\sqrt{2D}$. The scaling of this chi-squared error was chosen so that it has the same first two moments as the case with $\sigma = 6.4$ in Table 4(A). As expected, the results for large $D$ were essentially the same as those for $\sigma = 6.4$ in Table 4(A). Results for $D = 10$ are shown in the first two rows. The corresponding estimated powers at 1 and 5% of the general conditional coordinate test $T(\mathcal{H}|d=1)$ were observed to be 0.348 and 0.538. As illustrated in the third and fourth rows of this table, the performance of the marginal coordinate test is much better that the $t$-test when $D$ is small. The corresponding estimated powers at 1 and 5% of the

TABLE 4
*Power results based on model* (22) *with three different errors* $\varepsilon$

| Test | | $PR@0.01$ | $PR@0.05$ | $c_{0.05}$ | $c_{0.5}$ | $c_{0.95}$ |
|---|---|---|---|---|---|---|
| | | (A) $\varepsilon = \sigma N(0,1)$ | | | | |
| $t$ | $\sigma = 1$ | 1 | 1 | 0.977 | 0.992 | 0.998 |
| $T_n(\mathcal{H})$ | | 1 | 1 | 0.970 | 0.990 | 0.998 |
| $t$ | $\sigma = 6.4$ | 0.359 | 0.583 | 0.346 | 0.765 | 0.951 |
| $T_n(\mathcal{H})$ | | 0.175 | 0.364 | 0.095 | 0.583 | 0.772 |
| | | (B) $\varepsilon = 6.4(\chi^2(D) - D)/\sqrt{2D}$ | | | | |
| $t$ | $D = 10$ | 0.374 | 0.609 | 0.308 | 0.768 | 0.949 |
| $T_n(\mathcal{H})$ | | 0.220 | 0.465 | 0.120 | 0.698 | 0.948 |
| $t$ | $D = 2$ | 0.348 | 0.594 | 0.284 | 0.774 | 0.951 |
| $T_n(\mathcal{H})$ | | 0.797 | 0.928 | 0.605 | 0.895 | 0.976 |
| | | (C) $\varepsilon = (e^{\tau X_1})N(0,1)$ | | | | |
| $t$ | $\tau = 0.75$ | 1 | 1 | 0.959 | 0.987 | 0.997 |
| $T_n(\mathcal{H})$ | | 1 | 1 | 0.954 | 0.985 | 0.997 |
| $t$ | $\tau = 1.5$ | 0.508 | 0.630 | 0.177 | 0.817 | 0.980 |
| $T_n(\mathcal{H})$ | | 1 | 1 | 0.938 | 0.977 | 0.995 |

conditional coordinate test $T(\mathcal{H}|d=1)$ were observed to be 0.85 and 0.929. The results for the $t$- and $T_n(\mathcal{H})$-tests were found to be similar for $D$ around 5 or 6.

*Table* 4(C), $\varepsilon = (e^{\tau X_1})N(0,1)$. For $\tau$ near 0 this model is essentially the same as that for $\sigma = 1$ in Table 4(A), and the two tests were observed to be equivalent. However, with larger values of $\tau$ the $t$-test begins to lose ground and for sufficiently large values the performance of the coordinate test is again much better than the $t$-test. Results for $\tau = 0.75$ and 1.5 are shown.

The results of this section suggest that, while the coordinate tests might not perform as well as tests optimized for particular models, they perform reasonably across a wide range of regressions, particularly since they do not require a model for $Y|\mathbf{X}$.

7.3. *Choice of d*. As illustrated in the power study of Section 7.2, $T_n(\mathcal{H}|d)$ can be expected to have greater power than $T_n(\mathcal{H})$, and consequently there are potential gains from inferring about $d$ prior to testing predictors. On the other hand, misspecification of $d$ can lead to conclusions different from those based on the true value. In this section we describe qualitative results from a simulation study to investigate this behavior. Conclusions are based on the general marginal $T_n(\mathcal{H})$ and conditional $T_n(\mathcal{H}|d)$ coordinate test of each of the individual predictors.

Consider $n = 200$ observations from a regression with five independent standard normal predictors $X_j$ and response $Y = \mu(X_1, X_2) + \varepsilon$, where the standard normal error $\varepsilon \perp\!\!\!\perp \mathbf{X}$. When $\mu = X_1 + e^{X_2}$, the marginal dimension tests $T_n(m)$ resulted in the correct conclusion that $d = 2$, and the five conditional tests $T_n(\mathcal{H}|d=2)$ correctly concluded that only $X_1$ and $X_2$ are relevant to the regression. The five marginal tests $T_n(\mathcal{H})$ reached the same conclusion. With $d$ underspecified as 1, the five tests based on $T_n(\mathcal{H}|d=1)$ also resulted in the correct conclusion that only $X_1$ and $X_2$ are relevant. Underspecification did not affect the conclusions in this case because the first SIR direction was close to $\mathrm{Span}(\mathbf{e}_1 + \mathbf{e}_2)$, where $\mathbf{e}_i$ denotes the $5 \times 1$ vector with a 1 in the $i$th position and 0 otherwise. In other words, both $X_1$ and $X_2$ were manifested in the first SIR direction, and so $T_n(\mathcal{H}|d=1)$ was able to detect contributions from both predictors. With $d$ overspecified as 3, $T_n(\mathcal{H}|d=3)$ resulted in the conclusion that $X_1$, $X_2$ and $X_4$ are significant, thus giving an upper bound on the set of relevant predictors.

When $\mu = X_1 - X_2 + e^{(X_1+X_2)}$, the marginal dimension tests again resulted in the correct conclusion that $d = 2$, and $T_n(\mathcal{H}|d=2)$ and $T_n(\mathcal{H})$ again correctly concluded that only $X_1$ and $X_2$ are relevant to the regression. However, this time with $d$ underspecified as 1, the test $T_n(\mathcal{H}|d=1)$ incorrectly concluded that only $X_1$ is relevant. Underspecification affected the conclusions in this case because $X_2$ was not captured by the first SIR direction, which was close to $\mathrm{Span}(\mathbf{e}_1)$. With $d$ overspecified as 3, $T_n(\mathcal{H}|d=3)$ again indicated three significant predictors, including $X_1$ and $X_2$.

Results of this study, including results not reported here, suggest that misspecification of $d$ need not be a worrisome issue when the marginal dimension tests result in a clear estimate and that estimate is used in $T(\mathcal{H}|d)$. When the value of $d$ is not clear, it is still safe to base inference on the marginal coordinate test $T_n(\mathcal{H})$.

7.4. *Lean body mass regression.* We revisit the lean body mass regression [Cook and Weisberg (1999b)] to illustrate practical aspects of the previous development. Lean body mass (LBM) is regressed on the logarithms of height (Ht), weight (Wt), sum of skin folds (SSF) and the logarithms of the five hematological variables red cell count (RCC), white cell count (WCC), plasma ferritin concentration (PFC), hematocrit (Hc) and hemoglobin (Hg) for 202 athletes at the Australian Institute of Sport. Logarithms of the eight predictors were used to help insure the linearity condition. Both females and males are represented in the data in approximately equal proportions. However, for this illustration we neglect gender in the regression.

The SIR chi-squared $p$-values for the marginal dimension hypotheses $d = m, m = 0, 1, 2, 3,$ are about 0, 0, 0.13 and 0.46. Consequently, we initially inferred that $d = 2$, keeping in mind that $d = 3$ is also a possibility. The first two SIR directions $\hat{\boldsymbol{\eta}}_1$ and $\hat{\boldsymbol{\eta}}_2$ are shown in the second and third columns of Table 5. The numbers in parentheses are the approximate standard errors proposed by Chen and Li [(1998), page 297]. A scatterplot of LBM versus the first SIR predictor $\hat{\boldsymbol{\eta}}_1^T \mathbf{X}$ is shown in Figure 1(a). The mean function in this plot is noticeably curved. Letting $e$ denote the residuals from the OLS fit of LBM on $(\hat{\boldsymbol{\eta}}_1^T \mathbf{X}, \hat{\boldsymbol{\eta}}_2^T \mathbf{X})$, the need for a second direction is evident in a 3D plot of $e$ versus $(\hat{\boldsymbol{\eta}}_1^T \mathbf{X}, \hat{\boldsymbol{\eta}}_2^T \mathbf{X})$, which has a clear saddle shape. A scatterplot of $e$ versus $\hat{\boldsymbol{\eta}}_2^T \mathbf{X}$ is shown in Figure 1(b).

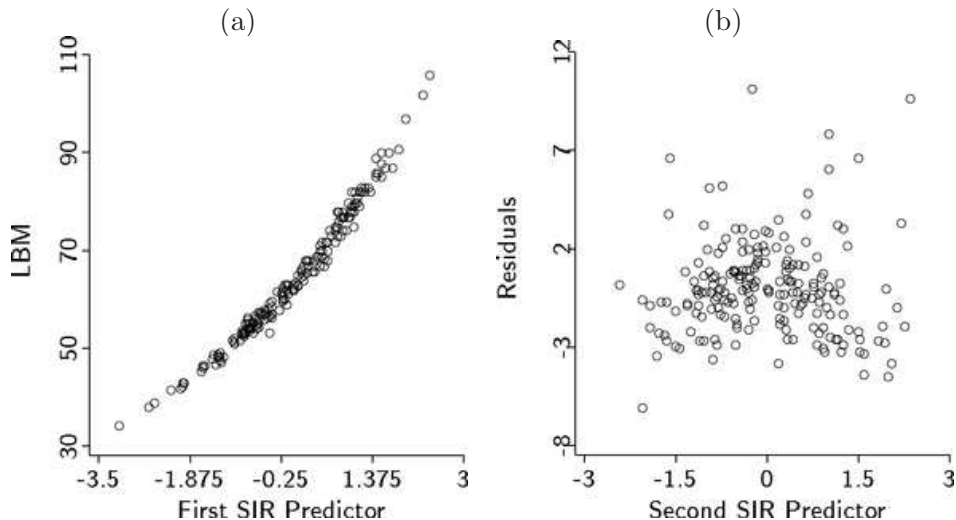

Fig. 1. *Two scatterplots representing the SIR "fit" of the lean body mass regression:* (a) *LBM versus* $\hat{\boldsymbol{\eta}}_1^T \mathbf{X}$*;* (b) *residuals versus* $\hat{\boldsymbol{\eta}}_2^T \mathbf{X}$.

TABLE 5
*Results from the lean body mass regression with all eight predictors*

| X | Fit | | $p$-values | | |
|---|---|---|---|---|---|
| | $\hat{\boldsymbol{\eta}}_1$ | $\hat{\boldsymbol{\eta}}_2$ | $T_n(\mathcal{H})$ | $T_n(\mathcal{H}|d=2)$ | $T_n(\mathcal{H}|d=3)$ |
| log[SSF] | −0.158 (0.06) | −0.076 (0.45) | 0 | 0 | 0 |
| log[Wt] | 0.971 (0.22) | −0.023 (1.6) | 0 | 0 | 0 |
| log[Hg] | 0.140 (0.69) | 0.347 (5.3) | 0.830 | 0.199 | 0.369 |
| log[Ht] | 0.088 (0.65) | −0.332 (5.0) | 0.344 | 0.270 | 0.537 |
| log[WCC] | −0.007 (0.08) | −0.015 (0.59) | 0.794 | 0.650 | 0.899 |
| log[RCC] | 0.011 (0.49) | 0.502 (3.8) | 0.090 | 0.014 | 0.032 |
| log[Hc] | −0.073 (0.85) | −0.715 (6.5) | 0.221 | 0.021 | 0.098 |
| log[PFC] | 0.003 (0.03) | 0.004 (0.25) | 0.040 | 0.820 | 0.192 |

In the context of SDR there are now at least three options to aid in assessing the significance of the individual predictors to the regression. We might develop a model for LBM|$\mathbf{X}$, guided by a 3D summary plot of LBM versus $(\hat{\boldsymbol{\eta}}_1^T \mathbf{X}, \hat{\boldsymbol{\eta}}_2^T \mathbf{X})$. Predictors could then be tested in the context of the resulting model. This type of procedure has produced useful results in the past, but there could be a worrisome possibility that the modeling process would effectively invalidate nominal characteristics of subsequent tests. Another possibility is to follow the case study by Chen and Li [(1998), Section 5.2] and use the approximate standard errors to guide variable selection. The assessment here is based on the general versions of the marginal and conditional coordinate tests.

The last three columns of Table 5 give the $p$-values from the marginal $T_n(\mathcal{H})$-test and the conditional tests $T_n(\mathcal{H}|d=2)$ and $T_n(\mathcal{H}|d=3)$ applied to each predictor in turn. We see from all three sets of tests that SSF and Wt contribute significantly to the regression, and probably RCC as well. The correlation between the first SIR predictor $\hat{\boldsymbol{\eta}}_1^T \mathbf{X}$ based on the full data and the first SIR predictor from the regression of LBM on log(SSF), log(Wt) is about 0.9995, so these two identified predictors largely account for the shape of the plot in Figure 1(a). The correlation between the second SIR predictors from the same regressions is about 0.83. Evidently, SSF and Wt contribute significantly to the first two directions, while other predictors contribute mostly to the second direction. As in linear regression, two correlated predictors might both have relatively large $p$-values, while deleting either causes the *p*-value for the remaining predictor to decrease substantially. Using $T_n(\mathcal{H})$ to test simultaneously the effects of the last six predictors in Table 5 yields a *p*-value of about 0.034, suggesting that some of those predictors also contribute to the regression. The tests $T_n(\mathcal{H}|d)$ with $d = 2, 3$ produced the same conclusion with similar *p*-values.

TABLE 6
*Results from the lean body mass regression with four predictors*

| | Fit | | | *p*-values | | |
|---|---|---|---|---|---|---|
| X | $\hat{\eta}_1$ | $\hat{\eta}_2$ | $\hat{\eta}_3$ | $T_n(\mathcal{H})$ | $T_n(\mathcal{H}|d=2)$ | $T_n(\mathcal{H}|d=3)$ |
| log[PFC] | 0.010 | 0.199 | 0.556 | 0.043 | 0.390 | 0.013 |
| log[RCC] | 0.023 | 0.556 | −0.714 | 0.004 | 0.039 | 0.001 |
| log[SSF] | −0.356 | −0.592 | −0.395 | 0 | 0 | 0 |
| log[Wt] | 0.934 | −0.549 | 0.159 | 0 | 0 | 0 |

The results so far can be partially summarized in terms of the hypothesis $Y \perp\!\!\!\perp \mathbf{X}_2|\mathbf{X}_1$, where $\mathbf{X}^T = (\mathbf{X}_1^T, \mathbf{X}_2^T)$. The tests gave firm indications that hypotheses with $\mathbf{X}_2 = \log(\text{SSF})$ and $\mathbf{X}_2 = \log(\text{Wt})$ are false. There was no notable information to reject hypotheses with $\mathbf{X}_2$ set to any combination of log(Ht), log(Hg) and log(WCC). Conclusions regarding the remaining three predictors were relatively ambiguous, depending on a dimension specification. One of the advantages of this type of analysis may be the ability to see which conclusions are firmly supported by the data without prespecifying a dimension for $\mathcal{S}_{Y|\mathbf{X}}$ and which depend on specification of a dimension and perhaps eventually a model. Nevertheless, to focus the analysis we deleted the three predictors that were judged to be unimportant and started over. The SIR chi-squared *p*-values from this regression for the hypotheses $d = m$, $m = 0, 1, 2, 3$, were about 0, 0, 0.010 and 0.31. Consequently, we now inferred that $d = 3$, a conclusion that remained stable for the rest of the analysis. This situation is consistent with the known propensity of the marginal dimension test to lose power when irrelevant predictors are added to the regression. Additionally, there was no notable evidence in this five-predictor regression to indicate that log(Hc) is relevant, leaving us with the reduced regression of LBM on the remaining four predictors (SSF, Wt, RCC, PFC).

Additional results for the reduced regression are given in Table 6. The sample correlations between the first, second and third SIR predictors from the full regression and the corresponding predictors from the reduced regression are 0.9997, 0.93 and 0.98, suggesting that the two regressions are giving essentially the same information about $Y|\mathbf{X}$. The *p*-values in Table 6 now give a consistent message for all predictors except PFC, which is judged nonsignificant when the dimension is underspecified as 2. This result could be anticipated from the discussion of underspecification in Section 7.3, if substantial information on PFC is furnished by the third direction. This interpretation is supported by the coefficients in Table 6, which were computed after marginally standardizing the predictors to have a sample standard deviation of 1. Additionally, the absolute sample correlations between log(PFC) and the three SIR predictors of Table 6 are about 0.32, 0.21 and 0.63.

In analogy with linear regression, we could have proceeded more straightforwardly by using backward elimination based on marginal or conditional tests to arrive at a reduced set of predictors. Starting with the marginal test in column 4, Table 5, and sequentially removing predictors whose $p$-values are larger than 0.05 yields the results given in the fifth column of Table 6. The same procedure based on the conditional test with $d=3$ yields the results in the last column of Table 6. The conditional test with $d=2$ ends in with the same predictors, except PFC is excluded because $d$ is underspecified.

**8. Discussion.** The theory of sufficient dimension reduction grew from a body of literature on how to graphically represent a regression in low dimensions without loss of information on $Y|\mathbf{X}$. Much of the development was inspired by the idea of linear sufficient statistics developed from a parametric view by Peters, Redner and Decell (1978) and Li's (1991) development of SIR. The primary motivation for the regression graphics ideas in Cook (1998a) stemmed from a desire to see how far graphics could be pushed in the analysis of regression data. The central subspace (CS) proved to be a key tool in that inquiry. The CS is intended to play a role similar to Li's (1991) EDR subspace, but it is a distinct population parameter constructed to insure that any nested sequence of dimension reduction subspaces always leads to the same population subspace. The methods developed here would not be possible using the EDR subspace because, in part, the fundamental equivalence of Proposition 1 would fail.

Data analytic techniques (e.g., SIR, SAVE, PIR) for pursuing sufficient dimension reduction have mostly lived in a world apart from mainstream methodology, although there are threads leading to other ideas and methods [see Chen and Li (1998) for a discussion]. By outlining a general context for testing predictors and developing a specific implementation using SIR, this article moves the inferential capabilities of SDR a step closer to mainstream regression methodology. The connection with tradition is also strengthened by casting SIR in terms of nonlinear least squares.

SIR has generated considerable interest since it was introduced. Hsing and Carroll (1992) develop a version of SIR in which each slice contains two observations so that the number of slices grows with the sample size. This two-slice method was extended by Zhu and Ng (1995) to allow for slices with more than two observations. The version in this article uses fixed slicing in which the number of observations per slice grows with the sample size. Zhu and Fang ([1996](#)) bypass the slicing step and use kernel smoothing instead. Schott ([1994](#)) investigated inference methods for $d$ when the predictors follow an elliptically contoured distribution. Elliptically contoured distributions are not required for the general methods in this article.

Cook and Critchley ([2000](#)) showed that SDR methods can be useful for identifying outliers and regression mixtures. Assuming $d$ to be known, Gather, Hilker and Becker ([2001](#)) developed a robust version of SIR by replacing its components (e.g., $\bar{\mathbf{X}}$, $\hat{\boldsymbol{\Sigma}}$ and $\widehat{\mathbf{M}}$) with robust estimates. The nonlinear least squares formulation of SIR described in Section [3](#) allows for alternative robust versions of SIR that involve using a loss function other than least squares.

The linearity condition (C1) and the coverage condition (C2) are the only two population conditions necessary for the theoretical justification of SIR. The constant covariance condition (C3) is used only to simplify the asymptotic distribution of the test statistic under the null hypothesis. Letting $\mathbf{M} = \mathrm{Var}(\mathrm{E}(\mathbf{Z}|Y))$, the role of the linearity and coverage conditions is to insure that $\mathrm{Span}(\mathbf{M}) = \mathcal{S}_{Y|\mathbf{Z}}$. Without these two conditions, the asymptotic distributions given in Theorems [1](#) and [2](#) remain valid if $\mathcal{S}_{Y|\mathbf{Z}}$ is replaced by $\mathrm{Span}(\mathbf{M})$, but we may lose the equality $\mathrm{Span}(\mathbf{M}) = \mathcal{S}_{Y|\mathbf{Z}}$ that provides an informative link with the population. As argued previously, the linearity condition need not be worrisome in practice, particularly if we use the adaptation methods discussed in Section [3.1](#). Li ([1997](#)) studied what can happen when the linearity condition fails, and Chen and Li ([1998](#)) developed an interpretation of SIR that might be helpful in some applications when $\mathrm{Span}(\mathbf{M}) \neq \mathcal{S}_{Y|\mathbf{Z}}$.

## APPENDIX: JUSTIFICATIONS

**A.1. Proposition [1](#).** Let the columns of the matrix $\boldsymbol{\eta}$ be a basis for $\mathcal{S}_{Y|\mathbf{X}}$. Then $Y \perp\!\!\!\perp \mathbf{X} | \boldsymbol{\eta}^T \mathbf{X}$ if and only if $Y \perp\!\!\!\perp (P_{\mathcal{H}}\mathbf{X}, Q_{\mathcal{H}}\mathbf{X}) | (\boldsymbol{\eta}^T P_{\mathcal{H}} \mathbf{X} + \boldsymbol{\eta}^T Q_{\mathcal{H}} \mathbf{X})$. Now,

$$P_{\mathcal{H}} \mathcal{S}_{Y|\mathbf{X}} = \mathcal{O}_p \quad \Longrightarrow \quad Y \perp\!\!\!\perp (P_{\mathcal{H}}\mathbf{X}, Q_{\mathcal{H}}\mathbf{X}) | \boldsymbol{\eta}^T Q_{\mathcal{H}} \mathbf{X}$$
$$\Longrightarrow \quad Y \perp\!\!\!\perp P_{\mathcal{H}}\mathbf{X} | Q_{\mathcal{H}}\mathbf{X}.$$

For the reverse implication, $Y \perp\!\!\!\perp P_{\mathcal{H}}\mathbf{X} | Q_{\mathcal{H}}\mathbf{X}$ implies $Y \perp\!\!\!\perp (P_{\mathcal{H}}\mathbf{X}, Q_{\mathcal{H}}\mathbf{X}) | Q_{\mathcal{H}}\mathbf{X}$ and consequently $\mathrm{Span}(Q_{\mathcal{H}})$ is a dimension reduction subspace. Since the central subspace is assumed to exist, any dimension reduction subspace must contain the central subspace, which therefore must be in $\mathrm{Span}(Q_{\mathcal{H}})$. It follows that $P_{\mathcal{H}} \mathcal{S}_{Y|\mathbf{X}} = \mathcal{O}_p$.

**A.2. Theorem [1](#).** The $y$th column $\sqrt{n}(\hat{\boldsymbol{\alpha}}^T \mathbb{Z}_n)_y$ of $\sqrt{n}\hat{\boldsymbol{\alpha}}^T \mathbb{Z}_n$ can be expressed as

$$\sqrt{n}(\hat{\boldsymbol{\alpha}}^T \mathbb{Z}_n)_y = \sqrt{n}(\boldsymbol{\alpha}_x^T \hat{\boldsymbol{\Sigma}}^{-1} \boldsymbol{\alpha}_x)^{-1/2} \boldsymbol{\alpha}_x^T \hat{g}_y \hat{\boldsymbol{\Sigma}}^{-1/2} \bar{\mathbf{Z}}_y.$$

Recalling that $J_y = 1$ if $Y$ is in slice $y$ and 0 otherwise, and that $\boldsymbol{\beta}_y = \boldsymbol{\Sigma}^{-1}\operatorname{Cov}(\mathbf{X}, J_y)$, it seems straightforward to verify that $\boldsymbol{\beta}_y = f_y \boldsymbol{\Sigma}^{-1/2} \operatorname{E}(\mathbf{Z}|Y = y)$ with OLS estimator

$$\hat{\boldsymbol{\beta}}_y = \hat{f}_y \hat{\boldsymbol{\Sigma}}^{-1/2} \bar{\mathbf{Z}}_y. \tag{24}$$

The linearity condition implies that $\boldsymbol{\beta}_y \in \mathcal{S}_{Y|\mathbf{X}}$ and thus under the coordinate hypothesis $\boldsymbol{\alpha}_x^T \boldsymbol{\beta}_y = 0$. Consequently, $\sqrt{n} \boldsymbol{\alpha}_x^T \hat{\boldsymbol{\beta}}_y$ converges in distribution and we have

$$\sqrt{n}(\hat{\boldsymbol{\alpha}}^T \mathbb{Z}_n)_y = \sqrt{n}(\boldsymbol{\alpha}_x^T \hat{\boldsymbol{\Sigma}}^{-1} \boldsymbol{\alpha}_x)^{-1/2} \hat{g}_y^{-1} \boldsymbol{\alpha}_x^T \hat{\boldsymbol{\beta}}_y$$
$$= \sqrt{n}(\boldsymbol{\alpha}_x^T \boldsymbol{\Sigma}^{-1} \boldsymbol{\alpha}_x)^{-1/2} g_y^{-1} \boldsymbol{\alpha}_x^T \hat{\boldsymbol{\beta}}_y + O_p(n^{-1/2}).$$

Li, Cook and Chiaromonte (2003) provided a general expansion for OLS estimators that is applicable to (24). Using this we get

$$\sqrt{n} \boldsymbol{\alpha}_x^T (\hat{\boldsymbol{\beta}}_y - \boldsymbol{\beta}_y) = n^{-1/2} \boldsymbol{\alpha}_x^T \boldsymbol{\Sigma}^{-1/2} \sum_{i=1}^n \mathbf{Z}_i \varepsilon_{yi} + O_p(n^{-1/2}),$$

where

$$\varepsilon_{yi} = J_{yi} - f_y - \boldsymbol{\beta}_y^T (\mathbf{X}_i - \operatorname{E}(\mathbf{X}))$$

is the population residual from the OLS regression of $J_y$ on $\mathbf{X}$. Thus, substituting we get

$$\sqrt{n}(\hat{\boldsymbol{\alpha}}^T \mathbb{Z}_n)_y = n^{-1/2} g_y^{-1} \boldsymbol{\alpha}^T \sum_{i=1}^n \mathbf{Z}_i \varepsilon_{yi} + O_p(n^{-1/2}),$$

where $\boldsymbol{\alpha} = \boldsymbol{\Sigma}^{-1/2} \boldsymbol{\alpha}_x (\boldsymbol{\alpha}_x^T \boldsymbol{\Sigma}^{-1} \boldsymbol{\alpha}_x)^{-1/2}$ as defined in (2).

Next define the $p \times h$ matrix

$$\mathbf{W}_n = \sum_{i=1}^n (\mathbf{Z}_i \varepsilon_{1i}, \ldots, \mathbf{Z}_i \varepsilon_{hi}) = \sum_{i=1}^n \mathbf{Z}_i \boldsymbol{\varepsilon}_i^T,$$

where $\boldsymbol{\varepsilon}_i$ is the $h \times 1$ vector with elements $\varepsilon_{yi}$, $y = 1, \ldots, h$. Define also the $r \times h$ matrix

$$\mathbf{V}_n = \boldsymbol{\alpha}^T \mathbf{W}_n \mathbf{D}_g^{-1}.$$

Then we have shown that

$$\sqrt{n} \operatorname{vec}(\hat{\boldsymbol{\alpha}}^T \mathbb{Z}_n) = n^{-1/2} \operatorname{vec}(\mathbf{V}_n) + O_p(n^{-1/2})$$
$$= (\mathbf{D}_g^{-1} \otimes \boldsymbol{\alpha}^T) n^{-1/2} \operatorname{vec}(\mathbf{W}_n) + O_p(n^{-1/2})$$
$$= (\mathbf{D}_g^{-1} \otimes \boldsymbol{\alpha}^T) n^{-1/2} \sum_{i=1}^n \boldsymbol{\varepsilon}_i \otimes \mathbf{Z}_i + O_p(n^{-1/2}).$$

Because $\boldsymbol{\varepsilon}$ contains OLS residuals, it is uncorrelated with any linear function of $\mathbf{X}$; in particular, $\mathrm{Cov}(\boldsymbol{\varepsilon}, \mathbf{Z}) = 0$. It follows that $\mathrm{E}(\boldsymbol{\varepsilon} \otimes \mathbf{Z}) = 0$ and therefore, by the multivariate central limit theorem, $n^{-1/2} \mathrm{vec}(\mathbf{W}_n)$ converges in distribution to a normal random vector with mean 0 and covariance matrix $\mathrm{Var}(\boldsymbol{\varepsilon} \otimes \mathbf{Z})$. Consequently, $n^{-1/2} \mathrm{vec}(\mathbf{V}_n)$ converges in distribution to a normal random vector with mean 0 and covariance matrix

$$\boldsymbol{\Omega}_{\mathcal{H}} = (\mathbf{D}_g^{-1} \otimes \boldsymbol{\alpha}^T) \mathrm{Var}(\boldsymbol{\varepsilon} \otimes \mathbf{Z})(\mathbf{D}_g^{-1} \otimes \boldsymbol{\alpha})$$
$$= \mathrm{E}(\mathbf{D}_g^{-1} \boldsymbol{\varepsilon} \boldsymbol{\varepsilon}^T \mathbf{D}_g^{-1} \otimes \boldsymbol{\alpha}^T \mathbf{Z} \mathbf{Z}^T \boldsymbol{\alpha}),$$

which is the desired conclusion.

**A.3. Corollary 1.**

A.4. *Equation* (13). Under the linearity condition, $\boldsymbol{\varepsilon}$ is a measurable function of $Y$ and $\boldsymbol{\Gamma}_1^T \mathbf{Z}$, and

$$\boldsymbol{\Omega}_{\mathcal{H}} = \mathrm{E}[\mathbf{D}_g^{-1} \boldsymbol{\varepsilon} \boldsymbol{\varepsilon}^T \mathbf{D}_g^{-1} \otimes \boldsymbol{\alpha}^T \mathrm{E}(\mathbf{Z}\mathbf{Z}^T|(Y, \boldsymbol{\Gamma}_1^T \mathbf{Z}))\boldsymbol{\alpha}].$$

The linearity and coverage conditions imply that $\mathrm{E}(\mathbf{Z}\mathbf{Z}^T|Y, \boldsymbol{\Gamma}_1^T \mathbf{Z}) = \mathrm{E}(\mathbf{Z}\mathbf{Z}^T|\boldsymbol{\Gamma}_1^T \mathbf{Z})$ and thus

$$\boldsymbol{\Omega}_{\mathcal{H}} = \mathrm{E}[\mathbf{D}_g^{-1} \boldsymbol{\varepsilon} \boldsymbol{\varepsilon}^T \mathbf{D}_g^{-1} \otimes \boldsymbol{\alpha}^T \mathrm{E}(\mathbf{Z}\mathbf{Z}^T|\boldsymbol{\Gamma}_1^T \mathbf{Z})\boldsymbol{\alpha}].$$

The linearity and constant covariance conditions imply that

$$\mathrm{E}(\mathbf{Z}\mathbf{Z}^T|\boldsymbol{\Gamma}_1^T \mathbf{Z}) = \mathrm{Var}(\mathbf{Z}|\boldsymbol{\Gamma}_1^T \mathbf{Z}) + \mathrm{E}(\mathbf{Z}|\boldsymbol{\Gamma}_1^T \mathbf{Z})\mathrm{E}(\mathbf{Z}|\boldsymbol{\Gamma}_1^T \mathbf{Z})^T$$
$$= Q_{\mathcal{S}_{Y|\mathbf{Z}}} + P_{\mathcal{S}_{Y|\mathbf{Z}}} \mathbf{Z}\mathbf{Z}^T P_{\mathcal{S}_{Y|\mathbf{Z}}}.$$

Consequently,

$$\boldsymbol{\alpha}^T \mathrm{E}(\mathbf{Z}\mathbf{Z}^T|\boldsymbol{\Gamma}_1^T \mathbf{Z})\boldsymbol{\alpha} = \boldsymbol{\alpha}^T Q_{\mathcal{S}_{Y|\mathbf{Z}}} \boldsymbol{\alpha} + \boldsymbol{\alpha}^T P_{\mathcal{S}_{Y|\mathbf{Z}}} \mathbf{Z}\mathbf{Z}^T P_{\mathcal{S}_{Y|\mathbf{Z}}} \boldsymbol{\alpha}$$
$$= \boldsymbol{\alpha}^T Q_{\mathcal{S}_{Y|\mathbf{Z}}} \boldsymbol{\alpha}$$
$$= \boldsymbol{\alpha}^T \boldsymbol{\alpha} = I_r,$$

where we have used the facts that under the coordinate hypothesis $P_{\mathcal{S}_{Y|\mathbf{Z}}} \boldsymbol{\alpha} = 0$ and $\mathcal{H} = \mathrm{Span}(\boldsymbol{\alpha}) \subseteq \mathrm{Span}(\boldsymbol{\Gamma}_0)$ (see Proposition 2). Thus, in summary to this point,

(25) $$\boldsymbol{\Omega}_{\mathcal{H}} = \mathrm{E}(\mathbf{D}_g^{-1} \boldsymbol{\varepsilon} \boldsymbol{\varepsilon}^T \mathbf{D}_g^{-1}) \otimes I_r.$$

To simplify $\mathrm{E}(\mathbf{D}_g^{-1} \boldsymbol{\varepsilon} \boldsymbol{\varepsilon}^T \mathbf{D}_g^{-1})$, let $\mathbf{J}$ and $\mathbf{f} = \mathrm{E}(\mathbf{J})$ denote the $h \times 1$ vectors with elements $J_y$ and $f_y$, and write the residual vector as

$$\boldsymbol{\varepsilon} = \mathbf{J} - \mathbf{f} - (\boldsymbol{\mu} \mathbf{D}_g)^T \mathbf{Z}.$$

In anticipation of expanding $E(\boldsymbol{\varepsilon}\boldsymbol{\varepsilon}^T)$ we have

$$E(\mathbf{J}-\mathbf{f})(\mathbf{J}-\mathbf{f})^T = E(\mathbf{JJ}^T) - \mathbf{ff}^T = \mathbf{D}_f - \mathbf{ff}^T,$$

$$E(\mathbf{Z}(\mathbf{J}-\mathbf{f})^T) = E(E(\mathbf{Z}|Y)\mathbf{J}^T) = \boldsymbol{\mu}\mathbf{D}_g.$$

Then

$$E(\boldsymbol{\varepsilon}\boldsymbol{\varepsilon}^T) = (\mathbf{D}_f - \mathbf{ff}^T) - 2\mathbf{D}_g\boldsymbol{\mu}^T\boldsymbol{\mu}\mathbf{D}_g + \mathbf{D}_g\boldsymbol{\mu}^T\boldsymbol{\mu}\mathbf{D}_g$$

$$= (\mathbf{D}_f - \mathbf{ff}^T) - \mathbf{D}_g\boldsymbol{\mu}^T\boldsymbol{\mu}\mathbf{D}_g,$$

$$E(\mathbf{D}_g^{-1}\boldsymbol{\varepsilon}\boldsymbol{\varepsilon}^T\mathbf{D}_g^{-1}) = I_r - \mathbf{gg}^T - \boldsymbol{\mu}^T\boldsymbol{\mu}$$

$$= Q_g - \boldsymbol{\mu}^T\boldsymbol{\mu}.$$

A.5. *Asymptotic distribution.* Since $\mathrm{Var}(\mathbf{Z}) = I_p$, we have

$$I_p = E(\mathrm{Var}(\mathbf{Z}|Y)) + \mathrm{Var}(E(\mathbf{Z}|Y))$$

$$= E(\mathrm{Var}(\mathbf{Z}|Y)) + \boldsymbol{\mu}\boldsymbol{\mu}^T$$

and consequently $I_p - \boldsymbol{\mu}\boldsymbol{\mu}^T \geq 0$, which implies that the eigenvalues of $\boldsymbol{\mu}\boldsymbol{\mu}^T$ are between 0 and 1. The nonzero eigenvalues of $\boldsymbol{\mu}\boldsymbol{\mu}^T$ are the same as the nonzero eigenvalues of $\boldsymbol{\mu}^T\boldsymbol{\mu}$ and thus $I_h - \boldsymbol{\mu}^T\boldsymbol{\mu} \geq 0$. Combining this with the identity $\boldsymbol{\mu}Q_g = \boldsymbol{\mu}$ we have $Q_g - \boldsymbol{\mu}^T\boldsymbol{\mu} = Q_g(I_h - \boldsymbol{\mu}^T\boldsymbol{\mu})Q_g \geq 0$. It follows that the eigenvalues of $Q_g - \boldsymbol{\mu}^T\boldsymbol{\mu}$ are nonnegative. The convergence in distribution follows immediately because $\delta_1 \geq \cdots \geq \delta_h = 0$, each with multiplicity $r$, are the eigenvalues of $\boldsymbol{\Omega}_{\mathcal{H}}$.

A.6. **Corollary 2.** It follows from Bura and Cook [(2001a), Theorem 2 and its justification; see also Cook (1998a), page 213] that, under the linearity, coverage and constant covariance conditions,

$$(\boldsymbol{\Psi}_0^T \otimes \boldsymbol{\Gamma}_0^T)\boldsymbol{\Delta}(\boldsymbol{\Psi}_0 \otimes \boldsymbol{\Gamma}_0) = \boldsymbol{\Psi}_0^T Q_g \boldsymbol{\Psi}_0 \otimes I_{p-d}.$$

Thus,

$$\boldsymbol{\Omega}_d' = \boldsymbol{\Psi}_0^T Q_g \boldsymbol{\Psi}_0 \otimes \mathbf{F}_{\mathcal{H}}.$$

The matrix $\mathbf{F}_{\mathcal{H}}$ is a symmetric idempotent matrix of rank $p - d - r$ [Proposition 2(iv)] and, from the discussion following Corollary 1, $\boldsymbol{\Psi}_0^T Q_g \boldsymbol{\Psi}_0$ is a symmetric idempotent matrix of rank $h - d - 1$. Consequently, $\boldsymbol{\Omega}_d'$ is a symmetric idempotent matrix of rank $(h - d - 1)(p - d - r)$ and the conclusion follows.

**A.7. Proposition 5.** To find the limiting distribution of $T_n(\mathcal{H}|d)$ under the coordinate hypothesis, we first use (19) and (16) to write

$$(26) \qquad T_n(\mathcal{H}|d) = T_n(\mathcal{H}) - n\|(I_{h-d} \otimes \mathbf{G}_\mathcal{H})\operatorname{vec}(\mathbf{U}_n)\|^2 + o_p(1)$$

and, because $\boldsymbol{\Psi}\boldsymbol{\Psi}^T = I_h$,

$$\sum_{y=1}^{h}\sum_{j=1}^{n_h}\|P_{\widehat{\mathcal{H}}}\bar{\mathbf{Z}}_y\|^2 = n\operatorname{trace}(P_{\widehat{\mathcal{H}}}\mathbb{Z}_n\boldsymbol{\Psi}\boldsymbol{\Psi}^T\mathbb{Z}_n^T P_{\widehat{\mathcal{H}}})$$

$$= n\operatorname{trace}(P_{\widehat{\mathcal{H}}}\mathbb{Z}_n\boldsymbol{\Psi}_1\boldsymbol{\Psi}_1^T\mathbb{Z}_n^T P_{\widehat{\mathcal{H}}})$$

$$+ n\operatorname{trace}(\boldsymbol{\Gamma}^T P_{\widehat{\mathcal{H}}}\mathbb{Z}_n\boldsymbol{\Psi}_0\boldsymbol{\Psi}_0^T\mathbb{Z}_n^T P_{\widehat{\mathcal{H}}}\boldsymbol{\Gamma}).$$

The second term in this expression for $T_n(\mathcal{H})$ can be represented using

$$\sqrt{n}\boldsymbol{\Gamma}^T P_{\widehat{\mathcal{H}}}\mathbb{Z}_n\boldsymbol{\Psi}_0 = \begin{pmatrix} \sqrt{n}\boldsymbol{\Gamma}_1^T P_{\widehat{\mathcal{H}}}\mathbb{Z}_n\boldsymbol{\Psi}_0 \\ \sqrt{n}\boldsymbol{\Gamma}_0^T P_{\widehat{\mathcal{H}}}\mathbb{Z}_n\boldsymbol{\Psi}_0 \end{pmatrix} = \begin{pmatrix} 0 \\ \sqrt{n}\mathbf{G}_\mathcal{H}\boldsymbol{\Gamma}_0^T\mathbb{Z}_n\boldsymbol{\Psi}_0 \end{pmatrix} + o_p(1).$$

The conclusion for the first coordinate relies on the fact that $P_{\widehat{\mathcal{H}}}\boldsymbol{\Gamma}_1$ converges to 0 in probability under the hypothesis and $\sqrt{n}\operatorname{vec}(\mathbb{Z}_n\boldsymbol{\Psi}_0)$ converges in distribution. The conclusion for the second coordinate follows by an argument similar to that used in Section 6.1. Recalling that $\mathbf{U}_n = \boldsymbol{\Gamma}_0^T\mathbb{Z}_n\boldsymbol{\Psi}_0$, we have

$$n\operatorname{trace}(\boldsymbol{\Gamma}^T P_{\widehat{\mathcal{H}}}\mathbb{Z}_n\boldsymbol{\Psi}_0\boldsymbol{\Psi}_0^T\mathbb{Z}_n^T P_{\widehat{\mathcal{H}}}\boldsymbol{\Gamma}) = n\operatorname{trace}(\mathbf{G}_\mathcal{H}\mathbf{U}_n\mathbf{U}_n^T\mathbf{G}_\mathcal{H}) + o_p(1)$$

$$= n\|(I_{(h-d)} \otimes \mathbf{G}_\mathcal{H})\operatorname{vec}(\mathbf{U}_n)\|^2 + o_p(1).$$

Combining this result with (26) we have

$$(27) \qquad \begin{aligned} T_n(\mathcal{H}|d) &= n\operatorname{trace}(P_{\widehat{\mathcal{H}}}\mathbb{Z}_n\boldsymbol{\Psi}_1\boldsymbol{\Psi}_1^T\mathbb{Z}_n^T P_{\widehat{\mathcal{H}}}) + o_p(1) \\ &= n\operatorname{trace}(\widehat{\boldsymbol{\alpha}}^T\mathbb{Z}_n\boldsymbol{\Psi}_1\boldsymbol{\Psi}_1^T\mathbb{Z}_n^T\widehat{\boldsymbol{\alpha}}) + o_p(1) \\ &= \|(\boldsymbol{\Psi}_1^T \otimes I_r)\sqrt{n}\operatorname{vec}(\widehat{\boldsymbol{\alpha}}^T\mathbb{Z}_n)\|^2 + o_p(1). \end{aligned}$$

**A.8. Corollary 4.** Under the linearity, coverage and constant covariance conditions, the covariance matrix of the asymptotic normal distribution of $(\boldsymbol{\Psi}_1^T \otimes I)\sqrt{n}\operatorname{vec}(\widehat{\boldsymbol{\alpha}}^T\mathbb{Z}_n)$ can be found by using Corollary 1:

$$\begin{aligned} \boldsymbol{\Omega}_{\mathcal{H}|d} &= (\boldsymbol{\Psi}_1^T \otimes I_r)\boldsymbol{\Omega}_\mathcal{H}(\boldsymbol{\Psi}_1 \otimes I_r) \\ &= (\boldsymbol{\Psi}_1^T \otimes I_r)((Q_g - \boldsymbol{\mu}^T\boldsymbol{\mu}) \otimes I_r)(\boldsymbol{\Psi}_1 \otimes I_r) \\ &= (\boldsymbol{\Psi}_1^T Q_g \boldsymbol{\Psi}_1 - \mathbf{D}_\lambda) \otimes I_r \\ &= (I_d - \mathbf{D}_\lambda) \otimes I_r, \end{aligned}$$

where $\mathbf{D}_\lambda = \mathbf{D}_s^2$ is the $d \times d$ matrix of the nonzero eigenvalues of $\boldsymbol{\mu}^T\boldsymbol{\mu}$. The third equality follows because

$$\boldsymbol{\mu}^T\boldsymbol{\mu} = \boldsymbol{\Psi}_1\mathbf{D}_s\boldsymbol{\Gamma}_1^T\boldsymbol{\Gamma}_1\mathbf{D}_s\boldsymbol{\Psi}_1^T$$

and thus

$$\mathbf{D}_\lambda = \boldsymbol{\Psi}_1^T \boldsymbol{\mu}^T \boldsymbol{\mu} \boldsymbol{\Psi}_1 = \mathbf{D}_s^2.$$

The final equality follows because $\boldsymbol{\mu} Q_g = \boldsymbol{\mu}$ and thus $\boldsymbol{\Gamma}_1 \mathbf{D}_s \boldsymbol{\Psi}_1^T Q_g = \boldsymbol{\Gamma}_1 \mathbf{D}_s \boldsymbol{\Psi}_1^T$, which implies that $\boldsymbol{\Psi}_1^T Q_g = \boldsymbol{\Psi}_1^T$.

The eigenvalues of $\boldsymbol{\Omega}_{\mathcal{H}|d}$ are $1 - \lambda_j$ with multiplicity $r$, $j = 1, \ldots, d$. Consequently,

$$T_n(\mathcal{H}|d) \xrightarrow{\mathcal{L}} \sum_{j=1}^{d} (1 - \lambda_j) \chi_j^2(r).$$

**Acknowledgments.** The author thanks Douglas Hawkins and the referees for helpful comments on earlier versions of this article.

SCHOOL OF STATISTICS
UNIVERSITY OF MINNESOTA
ST. PAUL, MINNESOTA 55108
USA
E-MAIL: dennis@stat.umn.edu